\documentclass[preprint,12pt]{elsarticle}
\usepackage{eepic}
\usepackage{amsmath,amssymb,amsthm}
\usepackage{multicol}

\usepackage{subfig}

\newtheorem{theorem}{Theorem}
\newtheorem{definition}{Definition}
\newtheorem{lemma}{Lemma}
\newtheorem{remark}{Remark}
\DeclareMathOperator{\LE}{LE}
\DeclareMathOperator{\LCE}{LCE}
\DeclareMathOperator{\LD}{LD}

\begin{document}
\begin{frontmatter}

 \title{The Lyapunov dimension formula for \\
 the global attractor of the Lorenz system}

 \author[spbu]{G.A.~Leonov}
 %\ead{leonov@math.spbu.ru}
 \author[spbu,fin]{N.V.~Kuznetsov\corref{cor}}
 %\ead{nkuznetsov239@gmail.com}
 \author[spbu]{N.A.~Korzhemanova}
 \author[spbu]{D.V.~Kusakin}
 %cortext[cor]{Corresponding author}

 \address[spbu]{Faculty of Mathematics and Mechanics,
 St. Petersburg State University, 198504 Peterhof,
 St. Petersburg, Russia}
 \address[fin]{Department of Mathematical Information Technology,
 University of Jyv\"{a}skyl\"{a}, \\ 40014 Jyv\"{a}skyl\"{a}, Finland}

 \begin{abstract}
 The exact Lyapunov dimension formula for the Lorenz system
 has been analytically obtained first due to G.A.Leonov in 2002
 under certain restrictions on parameters, permitting classical values.
 He used the construction technique of
 special Lyapunov-type functions developed by him in 1991 year.
 Later it was shown
 that the consideration of  larger class  of Lyapunov-type functions
 permits proving the validity of this formula
 for all parameters of the system such that
 all the equilibria of the system are hyperbolically unstable.
 In the present work it is proved the validity of the formula
 for Lyapunov dimension for a wider variety of parameters values,
 which include all parameters satisfying the classical physical limitations.
 One of the motivation of this work is the possibility of computing
 a chaotic attractor in the Lorenz system in the case of one unstable
 and two stable equilibria.
 \end{abstract}


\end{frontmatter}

\section{Introduction}
The exact Lyapunov dimension formula for the Lorenz system
 has been analytically obtained first due to G.A.Leonov in 2002 \cite{Leonov-2002}
 under certain restrictions on parameters, permitting classical values.
  In his work it was used the technique of special Lyapunov-type functions,
 which had been created in 1991 year \cite{Leonov-1991-Vest}
 and then was developed in \cite{LeonovB-1992,BoichenkoLR-2005}.
 Later in the works \cite{Leonov-2012-PMM,Leonov-2013-DAN-LD,LeonovPS-2013-PLA}
 it was shown
 that the consideration of a wider class of Lyapunov-type functions
 allows to provide the validity of the formula
 for such parameters of the Lorenz system
 that all its equilibria are hyperbolically unstable.

 In this study it is proved the validity of the formula
 under classical restrictions on the parameters.
 The motivation for this investigation is a numerical localization
 of chaotic attractor in the Lorenz system in the case of one unstable
 and two stable equilibria \cite{Sparrow-1982,YuC-2004}.

\section{The Lorenz system}
 Consider the classical Lorenz system suggested in the original work
 of Edward Lorenz \cite{Lorenz-1963}:
\begin{equation}\label{sys:Lorenz-classic}
 \left\{
 \begin{aligned}
 &\dot x= \sigma(y-x)\\
 &\dot y= r x-y-xz\\
 &\dot z=-b z+xy.
 \end{aligned}
 \right.
\end{equation}
E.~Lorenz obtained his system as a truncated model
of thermal convection in a fluid layer.
The parameters of this system are positive:
\[
 \sigma>0, \ \rho>0, \ b >0,
\]
because of their physical meaning
(e.g., $b = 4(1+a^2)^{-1}$ is positive and bounded).

Active study of the Lorenz system gave rise to
the appearance and subsequent consideration of various Lorenz-like systems
 (see, e.g., \cite{Celikovsky-1994,ChenU-1999,LuChen-2002,Tigan-2008}).
A recent discussion  of the equivalence of some Lorenz-like systems and
the possibility of universal consideration of their behavior
can be found, e.g. in \cite{LeonovK-2015-AMC,LeonovKKK-2015}.

Since the system is dissipative and generates a dynamical system
for $\forall t\geq 0$ (to verify this, it suffices to consider the Lyapunov function
$V(x,y,z) = \frac{1}{2}(x^2 + y^2 + (z - r - \sigma)^2)$;
see, e.g., \cite{Lorenz-1963,BoichenkoLR-2005}),
it possesses a global attractor
(a bounded closed invariant set, which is globally attractive)
\cite{Chueshov-2002-book,BoichenkoLR-2005}.

\begin{figure}[!ht]
\centering
\subfloat[]{
 \label{fig:1}
 \includegraphics[width=0.33\textwidth]{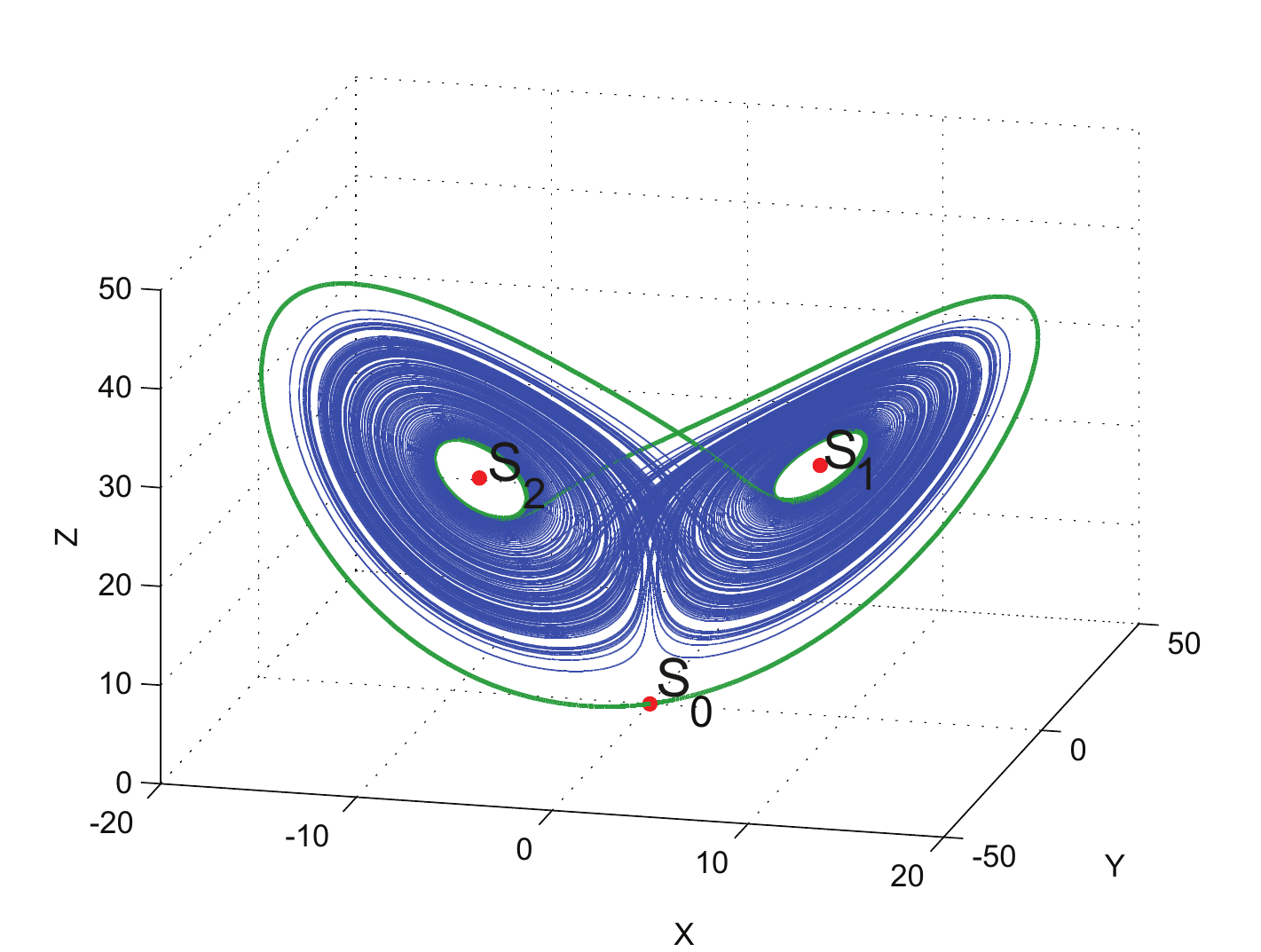}
}~
\subfloat[]{
 \label{fig:2}
 \includegraphics[width=0.33\textwidth]{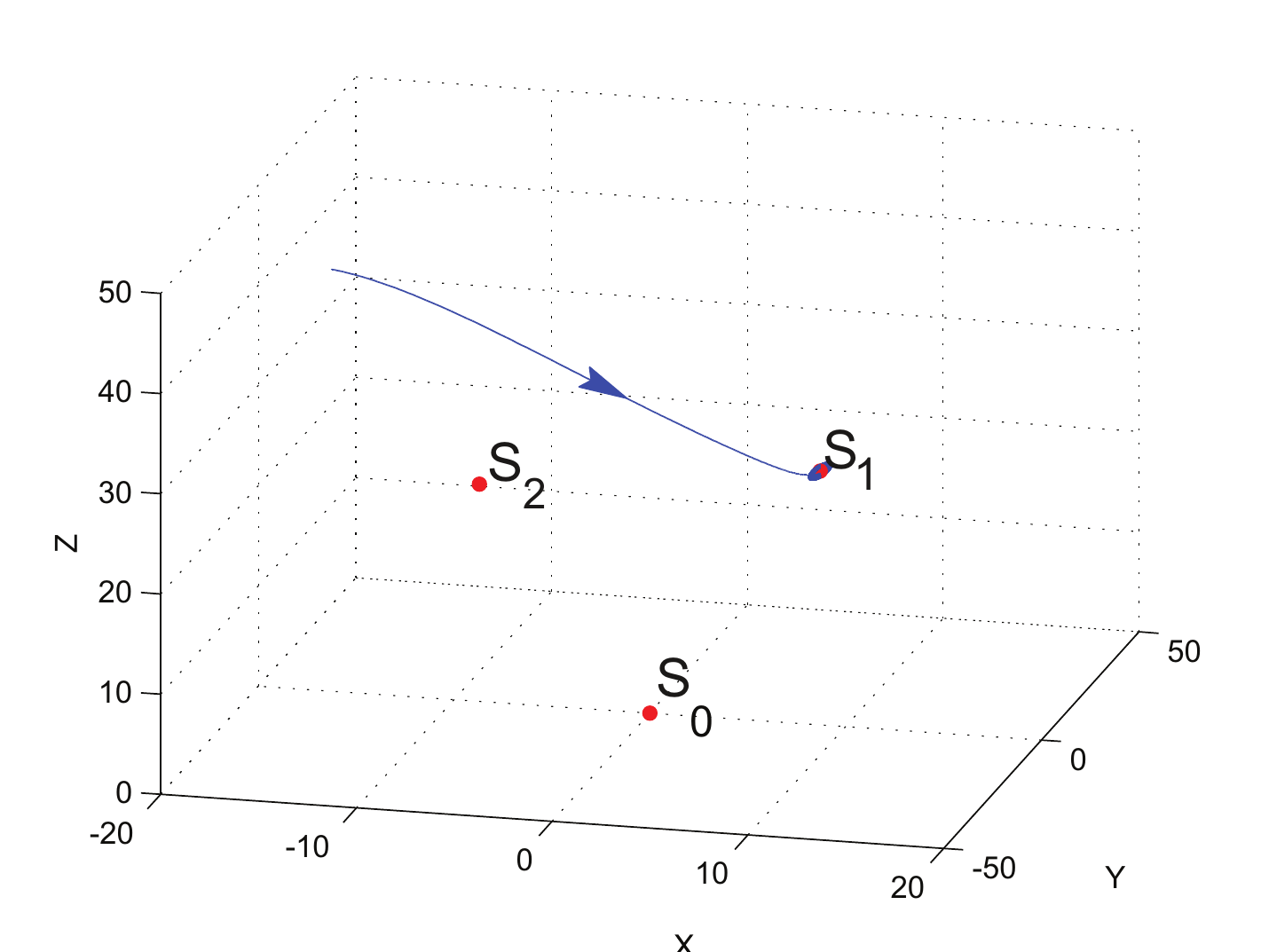}
}~
\subfloat[]{
 \label{fig:3}
 \includegraphics[width=0.33\textwidth]{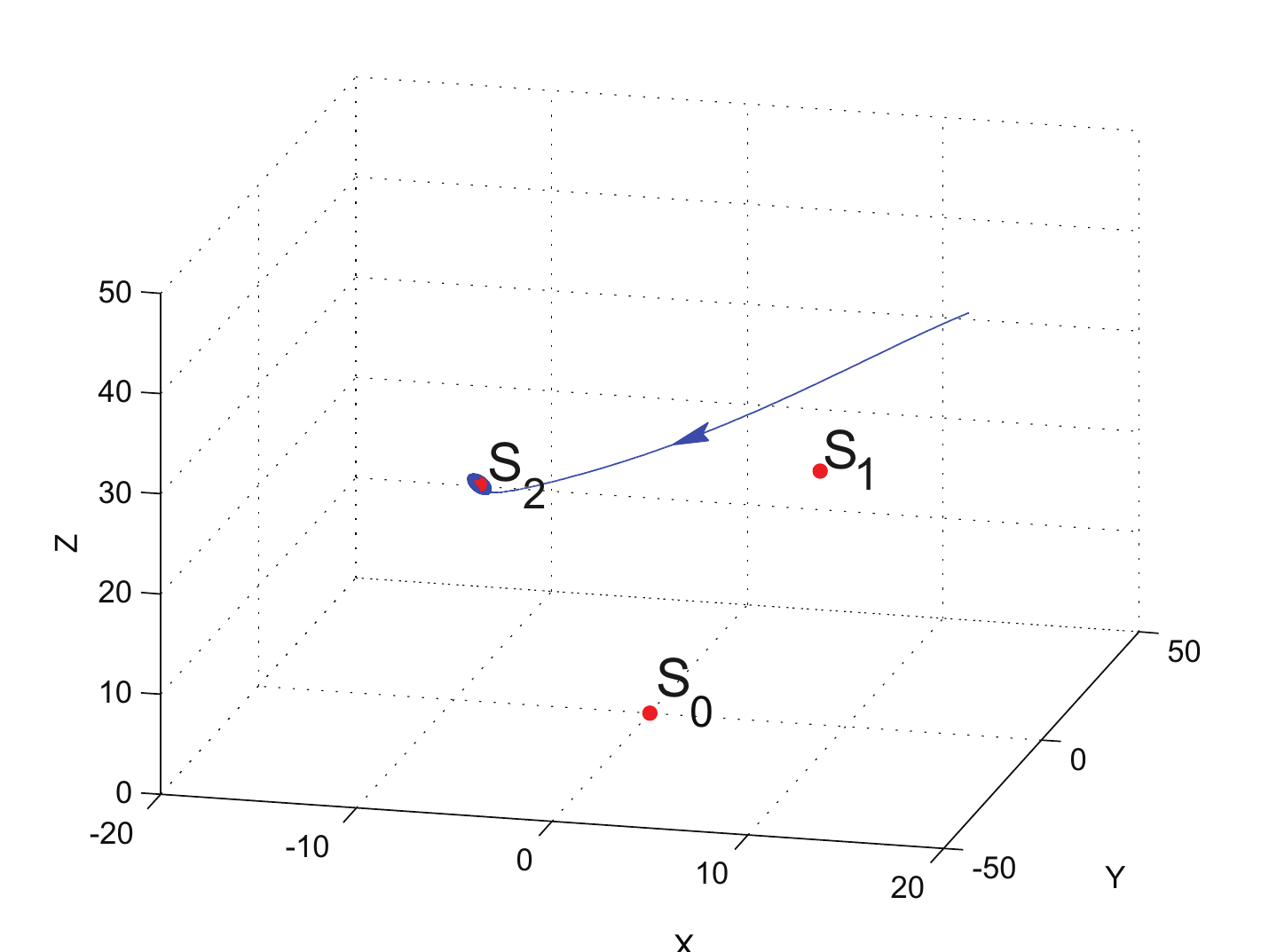}
}
\caption{
(a) Numerical visualization
of a self-excited Lorenz attractor
by using a trajectory with an initial point taken in
the vicinity of the equilibrium $S_0$.
(b), (c) Trajectories with the initial conditions
$(\mp 16.2899,\, \mp 0.0601,\, 42.1214)$ tend to
equilibria $S_{2,1}$.
Parameters: $r = 24.5$, $\sigma = 10$, $b = 8/3$.}
\label{fig:lorenz-co-exist}
\end{figure}

For the Lorenz system, the following classical scenario of transition
to chaos is known \cite{Sparrow-1982}.
Suppose that $\sigma$ and $b$ are fixed
(we use the classical parameters $\sigma = 10$, $b = 8/3$)
and $r$ varies.
Then, as $r$ increases,
the phase space of the Lorenz system is subject to the following sequence of bifurcations.
For $0<r<1$, there is globally asymptotically stable zero equilibrium $S_0$.
For $r>1$, equilibrium $S_0$ is a saddle
and a pair of symmetric equilibria $S_{1,2}$ appears.
For $1 < r < r_{h} \approx 13.9$, the separatrices $\Gamma_{1,2}$
of equilibria $S_0$ are attracted to the equilibria $S_{1,2}$.
For $r = r_{h} \approx 13.9$, the separatrices $\Gamma_{1,2}$
form two homoclinic trajectories of equilibria $S_0$ (homoclinic butterfly).
For $r_h < r < r_c \approx 24.06$, the separatrices $\Gamma_1$ and $\Gamma_2$
tend to $S_2$ and $S_1$, respectively.
For $r_c < r < r_a \approx 24.74$, the equilibria $S_{1,2}$ are stable and
the separatrices $\Gamma_{1,2}$ may be attracted to a local chaotic attractor
(see, e.g., \cite{Sparrow-1982,YuC-2004}).
This attractor is self-excited\footnote{
 An oscillation can generally be easily numerically localized if the initial data from its open
 neighborhood in the phase space (with the exception of a minor set of points) lead to
a long-term behavior that approaches the oscillation.
Therefore, from a computational perspective,
it is natural to suggest the following classification
of attractors \cite{KuznetsovLV-2010-IFAC,LeonovKV-2011-PLA,
LeonovKV-2012-PhysD,LeonovK-2013-IJBC},
which is based on the simplicity of finding their basins of attraction
in the phase space: {\emph{An attractor is called a self-excited attractor
 if its basin of attraction
 intersects with any open neighborhood of an equilibrium,
 otherwise it is called a hidden attractor}}
\cite{KuznetsovLV-2010-IFAC,LeonovKV-2011-PLA,LeonovKV-2012-PhysD,LeonovK-2013-IJBC}.
 Up to now in such Lorenz-like systems as
Lorenz, Chen, Lu and Tigan systems
 only self-excited chaotic attractors were found.
 In such Lorenz-like systems as Glukhovsky--Dolghansky and Rabinovich systems
 both self-excited and hidden attractors can be found
\cite{KuznetsovLM-2015,LeonovKM-2015-CNSNS,LeonovKM-2015-EPJST}.
Recent examples of hidden attractors
can be found in
\emph{The European Physical Journal Special Topics: Multistability: Uncovering Hidden Attractors}, 2015
(see \citep{Shahzad20151637,Brezetskyi20151459,Jafari20151469,Zhusubaliyev20151519,Saha20151563,Semenov20151553,Feng20151619,Li20151493,Feng20151593,Sprott20151409,Pham20151507,Vaidyanathan20151575}).
Note that while coexisting self-excited attractors
can be found by the standard computational procedure,
there is no regular way to predict the existence
or coexistence of hidden attractors.
%Recent examples of hidden attractors can be also found in
%\cite{
%KiselevaKLN-2012-IEEE,AndrievskyKLP-2013-IFAC,AndrievskyKLS-2013-IFAC,
%LeonovKKSZ-2014,KuznetsovLYY-2014-IFAC,KuznetsovKLNYY-2015-ISCAS,BestKKLYY-2015-ACC,
%SharmaSPKL-2015,
%DangLBW-2015-HA,
%KuznetsovKMS-2015-HA,
%ZhusubaliyevM-2015-HA,
%PhamVJWV-2014-HA,
%PhamJVWG-2014-HA,
%WeiWL-2014-HA,
%LiSprott-2014-HA,
%PhamRFF-2014-HA,
%WeiML-2014-HA,
%WeiYZY-2015-HA,
%ZhusubaliyevMRN-2015-HA,
%PhamVVLV-2015-HA,
%ChenYB-2015-HA,
%BaoHCXY-2015-HA,
%ChenLYBXW-2015-HA,
%WeiZWY-2015-HA,
%BurkinK-2014-HA,
%WeiZ-2014-HA,
%LiZY-2014-HA,
%ZhaoLD-2014-HA,
%LaoSJS-2014-HA,
%ChaudhuriP-2014-HA}
}
and can be found using the standard computational procedure, i.e.
by constructing a solution using initial data from a small neighborhood of zero equilibrium,
observing how it is attracted, and visualizing the attractor
(see Fig~\ref{fig:lorenz-co-exist}).

\begin{figure}[!ht]
 \centering
 \subfloat[
 {\scriptsize Initial data near the equilibrium $S_0$}
 ] {
 \label{fig:lorenz:attr:se0}
 \includegraphics[width=0.3\textwidth]{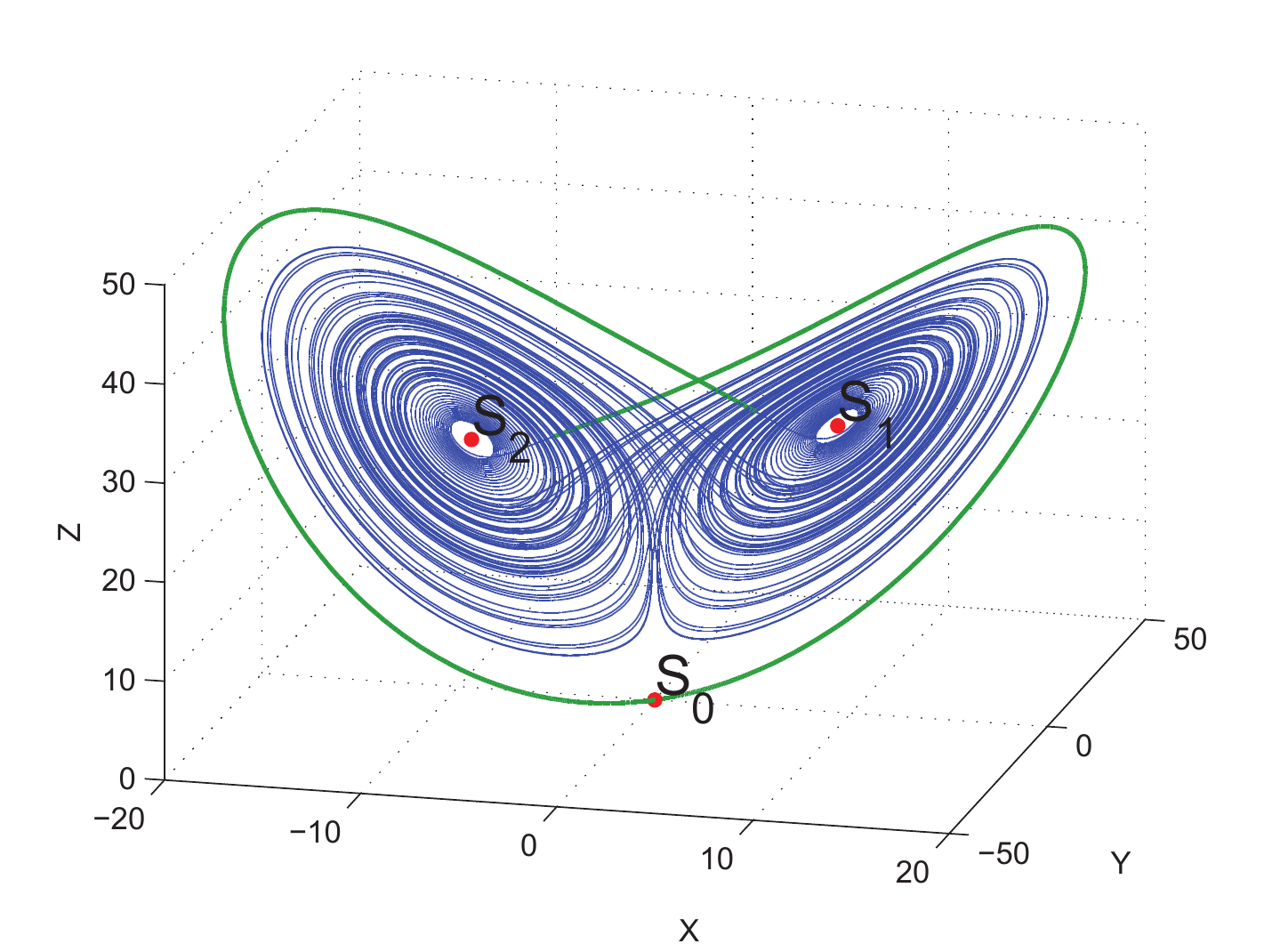}
 }~
 \subfloat[
 {\scriptsize Initial data near the equilibrium $S_1$}
 ] {
 \label{fig:lorenz:attr:se1}
 \includegraphics[width=0.3\textwidth]{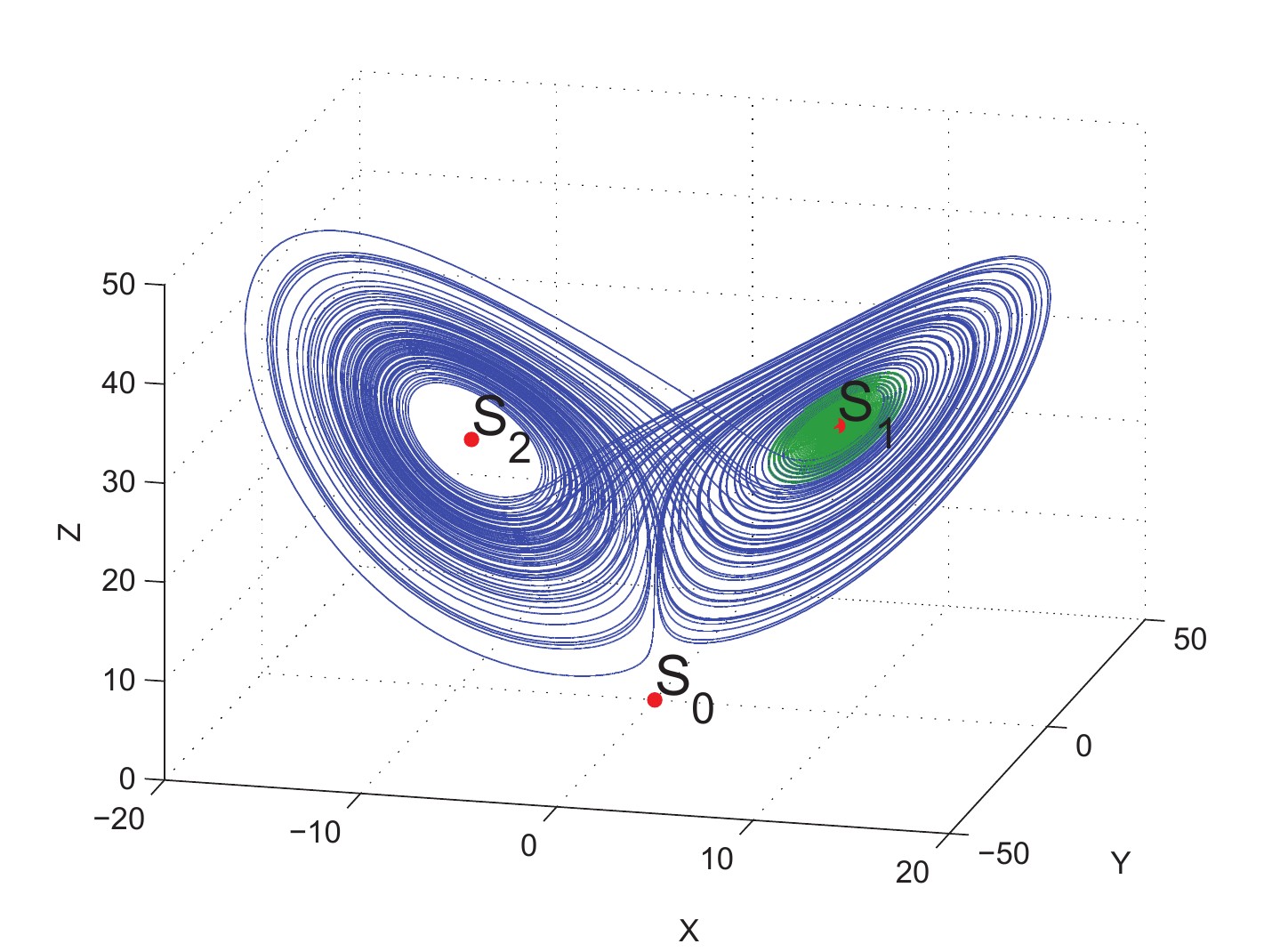}
 }~
 \subfloat[
 {\scriptsize Initial data near the equilibrium $S_2$}
 ] {
 \label{fig:lorenz:attr:se2}
 \includegraphics[width=0.3\textwidth]{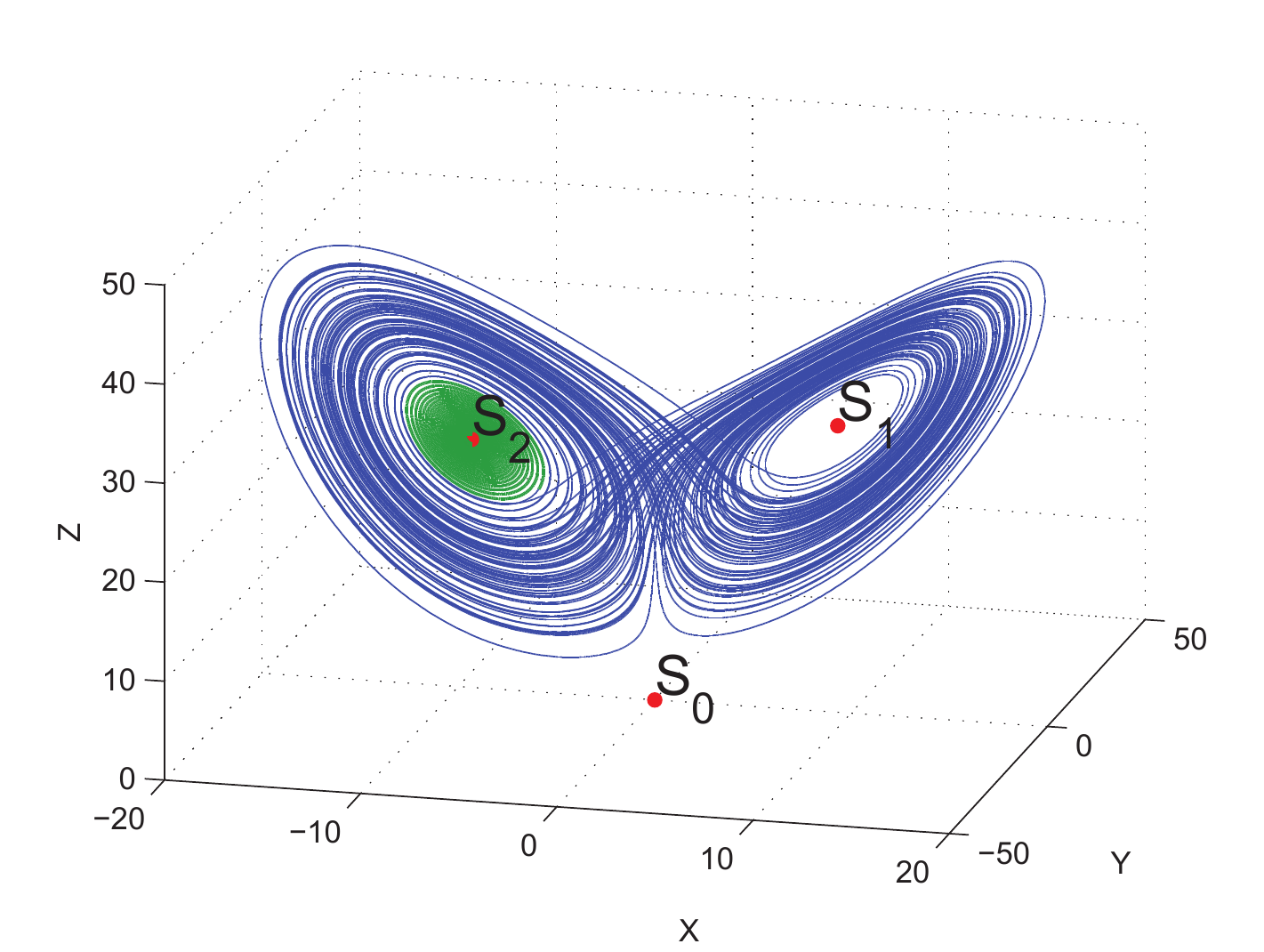}
 }
 \caption{
 Numerical visualization of
 the classical self-excited local attractor in the Lorenz system
 by using the trajectories that start in small neighborhoods of
the unstable equilibria $S_{0,1,2}$.
 Here the separation of the trajectory into transition process (green)
 and approximation of attractor (blue) is rough.
 }
 \label{fig:lorenz:attr:se}
\end{figure}

For $r>r_a$, the equilibria $S_{1,2}$ become unstable.
The value $r=28$ corresponds to the classical self-excited local attractor
(see Fig.~\ref{fig:lorenz:attr:se}).

\section{Lyapunov dimension of attractors}
Consider a dynamical system
\begin{equation}\label{eq:ode}
 \dot{{\bf \rm x}} = {\bf \rm f}({\bf \rm x}),
\end{equation}
where ${\bf \rm f} : \mathbb{R}^n \to \mathbb{R}^n$ is a smooth
vector-function. Since the initial time is not important for dynamical systems,
without loss of generality, we consider a solution
\(
 x(t,x_0):\ x(0,x_0)=x_0.
\)
The linearized system along the solution ${\bf \rm x}(t, {\bf \rm x}_0)$ is as follows
\begin{equation}\label{eq:lin_eq}
 \dot{{\bf \rm u}} = J({\bf \rm x}(t, {\bf \rm x}_0)) \, {\bf \rm u},
 \quad {\bf \rm u} \in \mathbb{R}^n, ~t \in \mathbb{R}_+,
\end{equation}
where
\[
J({\bf \rm x}(t, {\bf \rm x}_0)) =
\left[
\frac{\partial f_i({\bf \rm x})}{\partial {\bf \rm x}_j}\big|_{{\bf \rm x}
 = {\bf \rm x}(t, {\bf \rm x}_0)}
\right]
\]
is the $(n \times n)$ Jacobian matrix evaluated along the trajectory
${{\bf \rm x}(t, {\bf \rm x}_0)}$ of system \eqref{eq:ode}.
A fundamental matrix $X(t,{\bf \rm x}_0) $ of linearized system \eqref{eq:lin_eq}
is defined by the variational equation
\begin{equation}\label{eq:var_eq}
 \dot{X}(t,{\bf \rm x}_0) =
 J({\bf \rm x}(t, {\bf \rm x}_0)) \, X(t,{\bf \rm x}_0).
\end{equation}
We usually set $X(0,{\bf \rm x}_0) = I_n$, where $I_n$ is the
identity matrix.
Then
${\bf \rm u}(t, {\bf \rm u}_0) = X(t,{\bf \rm x}_0) {\bf \rm u}_0$.
In the general case,
${\bf \rm u}(t, {\bf \rm u}_0) = X(t,{\bf \rm x}_0)
X^{-1}(0,{\bf \rm x}_0){\bf \rm u}_0$.
Note that if a solution of nonlinear system \eqref{eq:ode}
is known, then we have
\[
X(t,{\bf \rm x}_0) =
\frac{\partial {\bf \rm x}(t, {\bf \rm x}_0)}{\partial {\bf \rm x}_0}.
\]

Consider the transformation of the unit ball $B$ into the ellipsoid $X(t,x_0)B$
and the exponential growth rates of its principal semiaxes lengths.
The principal semiaxes of the ellipsoid $X(t,x_0)B$  coincides with
singular values of the matrix $X(t,x_0)$,%:$\sigma_i(t,x_0) = \sigma_i(X(t,x_0))$,
which are defined as the square roots of the eigenvalues
of matrix $X(t,x_0)^*X(t,x_0)$.
Let $\sigma_1 (X(t,{\bf \rm x}_0)) \geq \cdots
\geq \sigma_n (X(t,{\bf \rm x}_0)) > 0$ denote
the singular values of the fundamental matrix $X(t,{\bf \rm x}_0)$
%(the square roots of the eigenvalues of the matrix
%$X(t,{\bf \rm x}_0)^{*}X(t,{\bf \rm x}_0)$ are
reordered for each $t$.
Introduce the operator $\mathcal{X}(\cdot)=\frac{1}{t}\ln|\cdot|$,
where $|\cdot|$ is the Euclidian norm.
Define the decreasing sequence (for all considered $t$)
of Lyapunov exponents functions $\LE_i(t,x_0) = \mathcal{X}({\sigma_i}(t,x_0))$
\[
 \LE_1(t,x_0)\geq \LE_2(t,x_0) \geq ... \geq \LE_n(t,x_0).
\]
\begin{definition}\label{def:le}\cite{Oseledec-1968}
 The {\it Lyapunov exponents (LEs)}
 at the point ${\bf \rm x}_0$
 are the numbers (or the symbols $\pm \infty$) defined as
 \begin{equation}\label{defLE}
 {\rm LE}_i ({\bf \rm x}_0) = \limsup_{t \to \infty} \LE_i(t,x_0)
 = \limsup_{t \to \infty} \frac{1}{t} \ln \sigma_i (X(t,{\bf \rm x}_0)).
 \end{equation}
\end{definition}

LEs are commonly used\footnote{
Two widely used definitions of Lyapunov exponents are
the upper bounds of the exponential growth rate
of the norms of linearized system solutions (LCEs) \cite{Lyapunov-1892}
and the upper bounds of the
exponential growth rate of the singular values
of fundamental matrix of linearized system (LEs) \cite{Oseledec-1968}.
 The LCEs \cite{Lyapunov-1892} and LEs \cite{Oseledec-1968} are ``often'' equal,
 e.g. for a ``typical'' system that satisfies the conditions
 of the Oseledec theorem \cite{Oseledec-1968}.
 However, there are no effective rigorous analytical methods for checking
 the Oseledec conditions for a given system \cite[p.118]{ChaosBook}
 (a numerical approach is discussed in \cite{DellnitzJ-2002}).
% :
% ``{\it Oseledec proof is important mathematics,
% but the method is not helpful in elucidating dynamics}'' \cite[p.118]{ChaosBook}).
 For particular system, LCEs and LEs may be different.
 For example,
 for the fundamental matrix
 \(
 X(t)=\left(
 \begin{array}{cc}
 1 & g(t)-g^{-1}(t) \\
 0 & 1 \\
 \end{array}
 \right)
 \) we have the following ordered values:
 $ \LCE_1 =
 {\rm max}\big(\limsup\limits_{t \to +\infty}\mathcal{X}[g(t)],
 \limsup\limits_{t \to +\infty}\mathcal{X}[g^{-1}(t)]\big),
 \LCE_2 = 0$;

 $
 \LE_{1,2} = {\rm max, min}
 \big(
 \limsup\limits_{t \to +\infty}\mathcal{X}[g(t)],
 \limsup\limits_{t \to +\infty}\mathcal{X}[g^{-1}(t)]
 \big)
 $. %, where $\mathcal{X}(\cdot) = \frac{1}{t}\log|\cdot|$
 %Thus, in general, Lyapunov dimensions, based
 %on LEs and LCEs, may be different.
 Note also that positive largest LCE or LE, computed via
 the linearization of the system along a trajectory,
 does not necessary imply instability or chaos,
 because for non-regular linearization
 there are well-known Perron effects of Lyapunov exponent sign reversal
 \cite{LeonovK-2007,KuznetsovL-2001,KuznetsovL-2005}.
 Therefore, in general, for the computation of the Lyapunov dimension
 of attractor we have to consider a grid of points on the attractor and
 corresponding local Lyapunov dimensions \cite{KuznetsovMV-2014-CNSNS}.
 More detailed discussion and examples can be found
 in \cite{KuznetsovAL-2014-arXiv-LE,LeonovK-2007}.
}
in the theory of dynamical systems and dimension theory
\cite{Ledrappier-1981,EckmannR-1985,Pesin-1988,Hunt-1996,Temam-1997,
BoichenkoLR-2005,BarreiraG-2011}.
\begin{remark}\label{exchangeLemmaLE}
The LEs are independent of the choice of fundamental matrix
at the point ${\bf \rm x}_0$ \cite{KuznetsovAL-2014-arXiv-LE,LeonovAK-2015}.
%unlike the Lyapunov characteristic exponents (LCEs, see \cite{Lyapunov-1892}).
%To determine all possible values of LCEs, we must consider a
%\emph{normal fundamental matrix}.
\end{remark}

Consider the largest integer $j(t,x_0) \in \{1,..,n\}$
%, which is equal to
%the largest natural number $m$
such that
\[
 \begin{array}{c}
 \LE_1^o(t,x_0) + \ldots + \LE_{j(t,x_0)}^o(t,x_0) > 0, \\ %\quad \LE_{j(t,x_0)+1}^o(t,x_0) < 0, \\
 \smallskip\\
 \cfrac{\LE_1^o(t,x_0) + \ldots + \LE_{j(t,x_0)}^o(t,x_0)}{|\LE_{j(t,x_0)+1}^o(t,x_0)|} < 1.
\end{array}
\]

Following \cite{KaplanY-1979,Hunt-1996},
introduce the following definition: the function $\LD(t,x_0) = 0$ if $\LE^{o}_1(t,x_0) \leq 0$
and $\LD(t,x_0)=n$ if $\sum_{i=1}^{n}\LE^{o}_i(t,x_0) \geq 0$,
otherwise
\begin{equation}\label{formula:kaplan}
 \LD(t,x_0) = j(t,x_0) + \cfrac{\LE_1^o(t,x_0) +
 \ldots + \LE_{j(t,x_0)}^o(t,x_0)}{|\LE_{j(t,x_0)+1}^o(t,x_0)|}.
\end{equation}
\begin{definition}
A local Lyapunov dimension at the point $x_0$ is as follows %\cite{KuznetsovAL-2014-arXiv-LE}
\[
 \LD(x_0) = \limsup\limits_{t \to +\infty} \LD(t,x_0).
\]
\end{definition}
The Lyapunov dimension of invariant compact set $K$ of
dynamical system
is defined by the relation
\[
  \dim_LK = \sup\limits_{x_0 \in K} \LD(x_0) = \sup\limits_{x_0 \in K} \limsup\limits_{t \to +\infty} \LD(t,x_0)
\].
%(for the case of exact limit this definition coincides
%with the classical one in \cite{KaplanY-1979}).

Note that, from an applications perspective,
an important property of the Lyapunov dimension is
the chain of inequalities~\cite{Hunt-1996,IlyashenkoW-1999,BoichenkoLR-2005}
\begin{equation}
 \dim_T K \leqslant \dim_H K \leqslant \dim_F K
 \leqslant \dim_L K \label{ineq:dimensions},
\end{equation}
where $\dim_T$, $\dim_H K,$ and $\dim_F K$ are
topological, Hausdorff, and fractal dimensions of $K$, respectively.

\section{Estimation of Lyapunov dimension by Lyapunov functions}
Along with commonly used numerical methods for estimating and computing
the Lyapunov dimension
(see, e.g., \cite{KuznetsovMV-2014-CNSNS,LeonovKM-2015-EPJST}),
there is an analytical approach that was proposed by G.A.Leonov
\cite{Leonov-1991-Vest,LeonovB-1992,BoichenkoLR-2005,Leonov-2008,
Leonov-2012-PMM,LeonovK-2015-AMC}.
It is based on the direct Lyapunov method
and uses Lyapunov-like functions.
The advantage of this method is that it allows
to estimate the Lyapunov dimension of invariant set
without numerical localization of the set.
This is especially important for the systems with hidden attractors
when numerical finding of all local attractors may be a challenging task
\cite{LeonovKM-2015-EPJST}.

Since LEs and LD are invariant under the linear
changes of variables (see, e.g., \cite{KuznetsovAL-2014-arXiv-LE}),
we can apply the linear variable change
${\bf\rm y} = S {\bf\rm x}$ with a nonsingular $n \times n$-matrix $S$.
Then system \eqref{eq:ode} is transformed into the system
\[
\dot{{\bf\rm y}} = S \,\dot{{\bf\rm x}} = S \,{\bf\rm f}
(S^{-1}{\bf\rm y}) = \tilde{{\bf\rm f}}({\bf\rm y}).
\]
Consider the linearization along the corresponding solution
${\bf\rm y}(t, {\bf\rm y}_0) = S {\bf\rm x}(t, S^{-1} {\bf\rm x}_0)$,
that is,
\begin{equation}\label{eq:lin_eq-new}
\dot {\bf \rm v} = \tilde{J}({\bf \rm y}(t, {\bf \rm y}_0))\,{\bf \rm v},
\quad {\bf \rm v} \in \mathbb{R}^n.
\end{equation}
Here the Jacobian matrix is as follows

\begin{align}
\tilde{J}({\bf\rm y}(t, {\bf\rm y}_0))
% &=
% \left[\frac{\partial \tilde{f}_i({\bf\rm y})}{\partial y_j}
% \big|_{{\bf\rm y} = {\bf\rm y}(t, {\bf\rm y}_0)} \right] =
% \left[\frac{\partial \left(S f_i(S^{-1}{\bf\rm y})\right)}{\partial y_j}
% \big|_{{\bf\rm y} = {\bf\rm y}(t, {\bf\rm y}_0)} \right] = \nonumber \\ &=
% S \left[ \frac{\partial f_i({\bf\rm x})}{\partial x_j}
% \big|_{{\bf\rm x} = {\bf\rm x}(t, {\bf\rm x}_0)} \right] S^{-1}
=S \, J({\bf\rm x}(t, {\bf\rm x}_0)) \, S^{-1} \label{jacobian-new}
\end{align}
and the corresponding fundamental matrix satisfies
\(
Y(t,{\bf \rm y}_0) = S X(t,{\bf \rm x}_0).
\)

For simplicity, let $J({\bf \rm x}) = J({\bf \rm x}(t, {\bf \rm x}_0))$.
Suppose that

$\lambda_1 ({\bf \rm x},S) \geqslant \cdots \geqslant \lambda_n ({\bf \rm x},S)$
are eigenvalues of the  symmetrized Jacobian matrix \eqref{jacobian-new}
\begin{equation}
 \frac{1}{2} \left( S J({\bf \rm x}) S^{-1} +
 (S J({\bf \rm x}) S^{-1})^{*}\right).
 \label{SJS}
\end{equation}

\begin{theorem}[\cite{Leonov-2012-PMM}]\label{theorem:th1}
 Given an integer $j \in [1,n]$ and $s \in [0,1]$,
 suppose that there are a continuously
 differentiable scalar function $\vartheta: \mathbb{R}^n \rightarrow \mathbb{R}$
 and a nonsingular matrix $S$ such that
 \begin{equation}\label{ineq:th-1}
 \lambda_1 ({\bf \rm x},S) + \cdots + \lambda_j ({\bf \rm x},S)
 + s\lambda_{j+1}
 ({\bf \rm x},S) + \dot{\vartheta}({\bf \rm x}) < 0,
 ~ \forall \, {\bf \rm x} \in K.
 \end{equation}
 Then $\dim_L K \leqslant j+s$.
\end{theorem}
Here $\dot{\vartheta}$ is the derivative of $\vartheta$ with respect
to the vector field ${\bf\rm f}$:
$$
\dot{\vartheta} ({\bf \rm x}) = ({\rm grad}(\vartheta))^{*}{\bf\rm f}({\bf \rm
 x}).
$$

\begin{theorem}[\cite{Leonov-1991-Vest,LeonovB-1992,BoichenkoLR-2005,Leonov-2012-PMM}]
 \label{theorem:th2}
 Assume that there are a continuously differentiable scalar function $\vartheta$
 and a nonsingular matrix $S$ such that
 \begin{equation}\label{ineq:th-2}
 \lambda_1 ({\bf \rm x},S) + \lambda_2 ({\bf \rm x},S) +
 \dot{\vartheta}({\bf \rm x}) < 0,
 ~ \forall \, {\bf \rm x} \in \mathbb{R}^n.
 \end{equation}
 Then any solution of system \eqref{eq:ode} bounded on $[0,+\infty)$
 tends to a certain equilibrium as $t \rightarrow +\infty$.
\end{theorem}

\section{Main result: Lyapunov dimension of the global Lorenz attractor}

 By Theorems \ref{theorem:th1} and \ref{theorem:th2},
for the Lorenz system we can obtain the following result.

\begin{theorem}\label{Main}
 Assume  that the following inequalities
 \begin{equation} \label{cond:mainTheoremR0}
 r - 1 > 0,
 \end{equation}
 \begin{equation}\label{cond:mainTheoremR}
 r - 1 \ge \frac{b(b + \sigma - 1)^2 -
 4\sigma(b + \sigma b - b^2)}{3\sigma^2}
 \end{equation}
 are satisfied.
 Let one of the following two conditions be satisfied:
 \begin{itemize}
 \item[a.]
 \begin{equation}\label{cond:mainTheorem1}
 \sigma^2 (r - 1)(b - 4)
 \le 4\sigma (\sigma b + b - b^2) - b(b + \sigma - 1)^2;
 \end{equation}

 \item[b.]
 there are two distinct real roots of equations
 \begin{equation}\label{eq:mainTheoremGamma2}
 \begin{split}
 (2\sigma - b + \gamma)^2 \left(b(b+\sigma-1)^2
 - 4\sigma(\sigma b + b-b^2) + \sigma^2(r-1)(b-4)\right) + \\
 + 4b\gamma(\sigma+1)\left(b(b+\sigma-1)^2
 - 4\sigma(\sigma b + b-b^2) - 3 \sigma^2(r-1)\right) =0
 \end{split}
 \end{equation}
 and
 \begin{equation}\label{cond:mainTheorem2}
 \left\{
 \begin{gathered}
 \sigma^2 (r - 1)(b - 4) > 4\sigma
 (\sigma b + b - b^2) - b(b + \sigma - 1)^2 \hfill \\
 \gamma^{(II)} > 0 \hfill \\
 \end{gathered}
 \right. \hfill\\
 \end{equation}
 where $\gamma^{(II)}$ is
 a greater root of equation \eqref{eq:mainTheoremGamma2}.
 \end{itemize}

 In this case
 \begin{enumerate}
 \item if
 \begin{equation}\label{cond:mainTheoremR1}
 (b-\sigma)(b-1)< \sigma r < (b + 1)(b + \sigma),
 \end{equation}
 then any bounded on $[0; +\infty)$
 solution of system \eqref{sys:Lorenz-classic}
 tends to a certain equilibrium as $t \to +\infty$.

 \item If
 \begin{equation}\label{cond:mainTheoremR2}
 \sigma r > (b + 1)(b+\sigma),
 \end{equation}
 then
 \begin{equation}\label{formula}
 \dim_L K \le 3 - \frac{2 (\sigma + b + 1)}{\sigma + 1
 + \sqrt{(\sigma-1)^2 + 4 \sigma r}},
 \end{equation}
 \end{enumerate}
 where $K$ is a bounded invariant set.
\end{theorem}

%\section{Computation of Lyapunov dimension at the point (0,0,0)}

For numerical experiments (see, e.g., \cite{EdenFT-1991,DoeringG-1995})
it is known that the Lyapunov dimension of any invariant compact set of the Lorenz system
is bounded from above by the local Lyapunov dimension of the zero equilibrium.

\begin{lemma}
 If for the parameters of system \eqref{sys:Lorenz-classic}
 inequality \eqref{cond:mainTheoremR2} is valid, then
 $$\dim_L(0,0,0)=3 - \frac{2(\sigma +b+1)}{\sigma+1
 + \sqrt{(\sigma-1 )^2 + 4r\sigma}}.$$
\end{lemma}

Comparing estimation \eqref{formula} with the local dimension
in the origin, we obtain a dimension formula
for a global Lorenz attractor.

\begin{theorem}\label{thm4}
 Let $K$ be a global attractor of classical Lorenz system
 \eqref{sys:Lorenz-classic}
 and conditions \eqref{cond:mainTheoremR0}-\eqref{cond:mainTheorem2}
 of Theorem \ref{Main} be hold true.

 In this case if
 \begin{equation*}
 \sigma r > (b + 1)(b+\sigma),
 \end{equation*}
 then
 \begin{equation}\label{mainformula}
 \dim_L K = 3 - \frac{2 (\sigma + b + 1)}{\sigma + 1 +
 \sqrt{(\sigma-1)^2 + 4 \sigma r}}.
 \end{equation}
\end{theorem}

\begin{remark}\label{r:unstableEquilibria}
 It can be easily checked numerically that if all three equilibria are
 hyperbolic, then the conditions of Theorem \ref{thm4} (see Fig.~\ref{fig:bounded-domain-r-28})
 are satisfied.
 For example, for the standard parameters $\sigma = 10$ and $b=\frac{8}{3}$
 the formula \eqref{mainformula} is valid for $r > \frac{209}{45}$.
 \begin{figure}[!ht]
 \centering
 \subfloat[
 {\scriptsize $r = 10$}
 ] {
 \label{fig:bounded-domain-r-10}
 \includegraphics[angle=270, width=0.33\columnwidth]{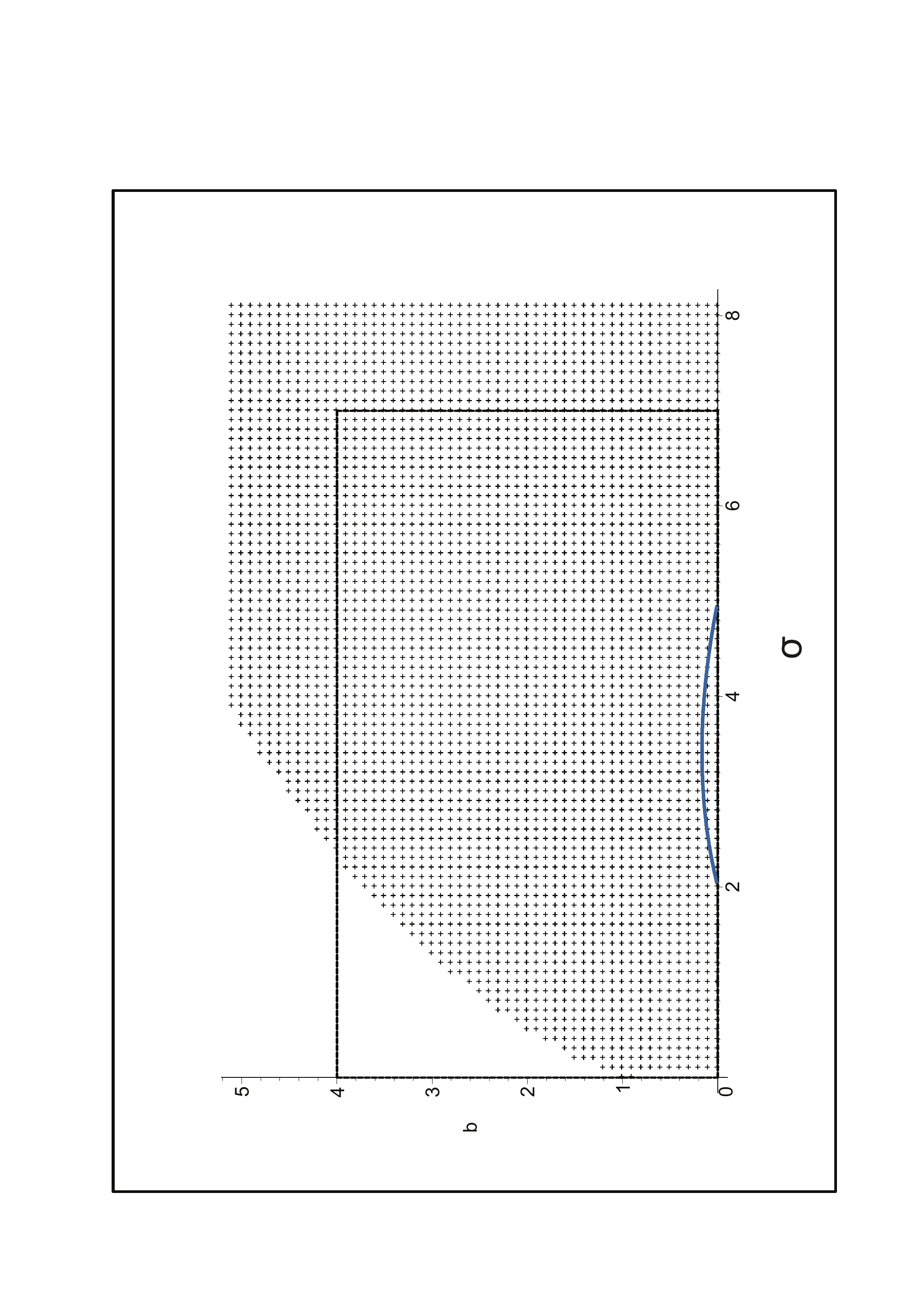}
 }~
 \subfloat[
 {\scriptsize $r = 28$}
 ] {
 \label{fig:bounded-domain-r-28}
 \includegraphics[angle=270, width=0.33\columnwidth]{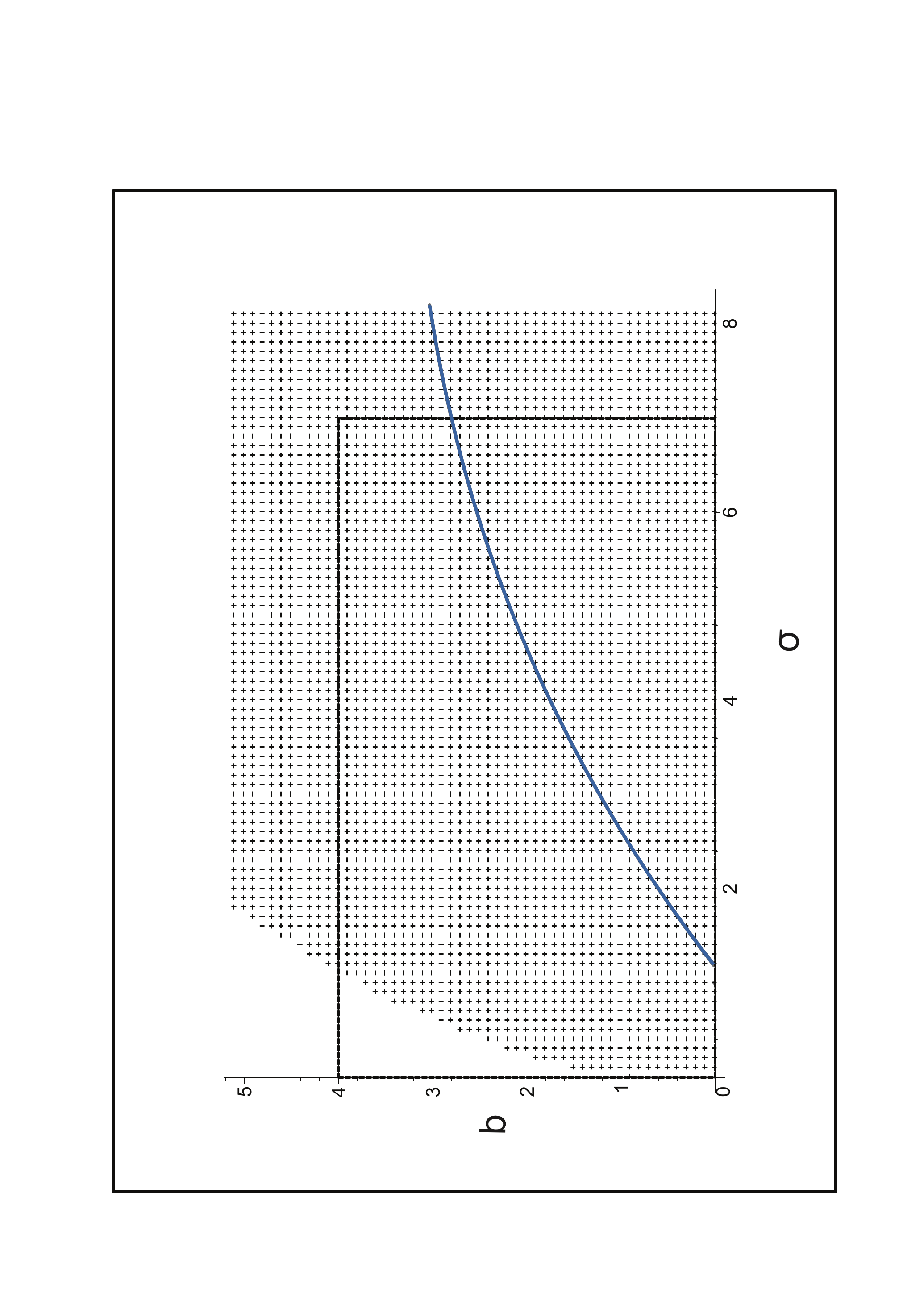}
 }~
 \subfloat[
 {\scriptsize $r = 100$}
 ] {
 \label{fig:bounded-domain-r-100}
 \includegraphics[angle=270, width=0.33\columnwidth]{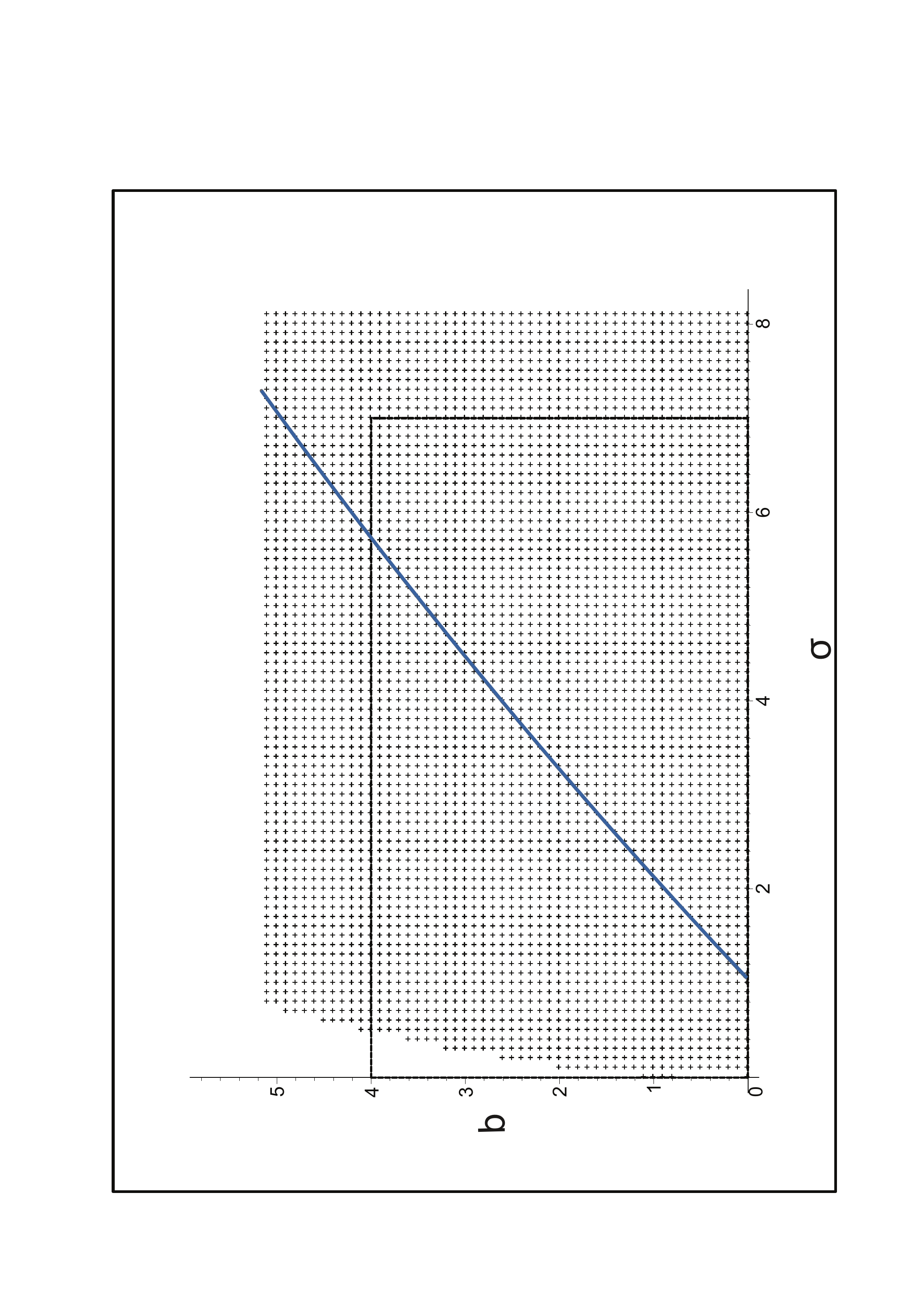}
 }
 \caption{
 Numerical simulation of domains for the fixed values of parameter
 $r= 10, 28, 100$
 }
 \label{fig:bounded-domains}
 \end{figure}
\end{remark}

\subsection{Numerical analysis}

We perform a numerical analysis of  parameter domain, which does not satisfy
the conditions of Theorem \ref{thm4}.

{\bf 1.} Recall that for $0<r<1$
the solutions of system tend to a unique equilibria
(see., e.g., \cite{BoichenkoLR-2005}).

{\bf 2.} Consider $r>1$.
 In the case of conditions \eqref{cond:mainTheoremR0}-\eqref{cond:mainTheoremR2}
 the following assertion is valid.
\begin{lemma}
 For all $r>1$, $\sigma>7$, and $0<b<4$
 the conditions of Theorem \ref{Main} are satisfied.
\end{lemma}

\begin{remark}
 For $0 < b < 4 - \epsilon$, $\sigma>7$ and $r > r_*$, where $r_* = r_*(\epsilon)$,
 the conditions of Theorem \ref{Main} are satisfied.

\end{remark}

Thus, for each fixed $r>1$ it remains to check numerically
the behavior of system with the parameters from the bounded
domain $[0<b<4, 0<\sigma<7]$. For the considered increasing sequence $r_k>1$,
numerical simulation shows the lack of chaotic attractor
(the trajectories with initial data from the corresponding compact
absorbing set \cite{LeonovBK-1987} have been simulated).

\section*{Appendix}

\begin{proof} (Theorem 1)

For the proof we make use of Theorems \ref{theorem:th1},
\ref{theorem:th2}.
Consider a solution of system \eqref{sys:Lorenz-classic}
${\bf \rm x}=(x(t), y(t), z(t))=(x,y,z)$.

\textbf{I.}
In the case of classical Lorenz system the matrix $J$ has the following form
$$J = \begin{pmatrix}
 -\sigma & \sigma & 0\\
 r - z & -1 & -x\\
 y & x & -b
\end{pmatrix}
$$

Following \cite{BoichenkoLR-2005}, for the condition
\begin{equation}\label{ineq:sqrt}
 \sigma r + (\sigma-b)(b-1) > 0
\end{equation}
we introduce a matrix
$$S = \begin{pmatrix}
 -\rho^{-1} & 0 & 0\\
 -\frac{b-1}{\sigma} & 1 & 0\\
 0 & 0 & 1
\end{pmatrix},$$
where
\begin{equation}\label{rho}
 \rho = \frac{\sigma}{\sqrt{\sigma r + (\sigma-b)(b-1)}}.
\end{equation}

Then
$$S J S^{-1}=\begin{pmatrix}
 b-\sigma-1 & -\frac{\sigma}{\rho} & 0\\
 -\frac{\sigma}{\rho} + \rho z & -b & -x\\
 -\rho \left(y + \frac{b-1}{\sigma} x\right) & x & -b
\end{pmatrix}.$$
We find the eigenvalues of the matrix
$\frac{1}{2} (S J S^{-1} + (S J S^{-1})^{*})$.
The characteristic polynomial of the matrix takes the form:
$$(\lambda + b) \left[\lambda ^ 2 + (\sigma+1) \lambda + b (\sigma - b + 1)
- \left(\frac{\sigma}{\rho} - \frac{\rho z}{2}\right)^2
-\left(\frac{\rho (b - 1)}{2 \sigma}x + \frac{\rho}{2}y\right)^2\right].$$

This implies that the eigenvalues $\lambda_i=\lambda_i(x,y,z,S), i=1,2,3$
of the matrix $\frac{1}{2} (S J S^{-1} + (S J S^{-1})^{*})$ are the numbers:
$$\lambda_2 = -b,$$
$$\lambda_{1,3} = -\frac{\sigma+1}{2} \pm \frac{1}{2}\left[(\sigma - 2 b
+ 1)^2 + \left(\frac{2 \sigma}{\rho}-\rho z\right)^2 + \rho^2
\left(y + \frac{b - 1}{\sigma} x\right)^2\right]^{\frac{1}{2}}.$$

Such a choice of the matrix $S$ provides a simple form of eigenvalues.

\bigskip
\textbf{II.}
We shall show that $\lambda_1 \ge \lambda_2 \ge \lambda_3, \forall x,y,z$:
\begin{align*}
 &\frac{-(\sigma+1) + \sqrt{(\sigma-2b+1)^2 + \left(\frac{2\sigma}{\rho}
 - \rho z \right)^2 + \rho^2 \left(y + \frac{b-1}{\sigma}x \right)^2}}{2} \ge -b \ge\\
 & \ge \frac{-(\sigma+1) - \sqrt{(\sigma-2b+1)^2 + \left(\frac{2\sigma}{\rho}
 -\rho z \right)^2 +\rho^2 \left(y + \frac{b-1}{\sigma}x \right)^2}}{2}
\end{align*}
%$$\Leftrightarrow$$
%\begin{align*}
% &-(\sigma+1) + \sqrt{(\sigma-2b+1)^2 + \left(\frac{2\sigma}{\rho} - \rho z \right)^2 + \rho^2 \left(y + \frac{b-1}{\sigma}x \right)^2} \ge -2b \ge \\
% & \ge -(\sigma+1) - \sqrt{(\sigma-2b+1)^2 + \left(\frac{2\sigma}{\rho} -\rho z \right)^2 +\rho^2 \left(y + \frac{b-1}{\sigma}x \right)^2}
%\end{align*}
$$\Leftrightarrow$$
\begin{align*}
 &\sqrt{(\sigma-2b+1)^2 + \left(\frac{2\sigma}{\rho} - \rho z \right)^2 + \rho^2 \left(y + \frac{b-1}{\sigma}x \right)^2} \ge \sigma-2b+1 \ge \\
 &\ge - \sqrt{(\sigma-2b+1)^2 + \left(\frac{2\sigma}{\rho} -\rho z \right)^2 +\rho^2 \left(y + \frac{b-1}{\sigma}x \right)^2}
\end{align*}
$$\Leftrightarrow$$
$$\sqrt{(\sigma-2b+1)^2 + \left(\frac{2\sigma}{\rho} -\rho z \right)^2 +\rho^2 \left(y + \frac{b-1}{\sigma}x \right)^2} \ge |\sigma-2b+1|$$
%both ðart@s @ side@s inequality@ies are positive, consequently their can @ it is possible âîçâåñòè squaring:
%$$(\sigma-2b+1)^2 + \left(\frac{2\sigma}{\rho} -\rho z \right)^2 +\rho^2
%\left(y + \frac{b-1}{\sigma}x \right)^2 \ge (\sigma-2b+1)^2$$
%$$\Leftrightarrow$$
$$\left(\frac{2\sigma}{\rho} -\rho z \right)^2 +\rho^2
\left(y + \frac{b-1}{\sigma}x \right)^2 \ge 0.$$

At the left there is a sum of two total squares, i.e. the inequality
is always satisfied.

Thus, $\lambda_1 \ge \lambda_2 \ge \lambda_3, \forall x,y,z$.

\bigskip
\textbf{III.}
To apply Theorem \ref{theorem:th1} and Theorem \ref{theorem:th2},
for $s \in [0,1)$, we consider and transform the relation
$2 (\lambda_1+\lambda_2+s\lambda_3)$:
\begin{equation}\label{transformation}
 \begin{split}
&2 (\lambda_1+\lambda_2+s\lambda_3) = -(\sigma+2b+1) - s (\sigma+1) + \\
&+(1-s)\left[(\sigma-2b+1)^2 + \left(\frac{2\sigma}{\rho}
-\rho z \right)^2 +\rho^2 \left(y + \frac{b-1}{\sigma}x \right)^2\right]^{\frac{1}{2}}=\\
&-(\sigma+2b+1) - s (\sigma+1) + \\
&+(1-s)\left[(\sigma-2b+1)^2+\frac{4\sigma^2}{\rho^2}-4\sigma z+\rho^2 z^2
+ \rho^2 \left(y + \frac{b-1}{\sigma}x\right)^2\right]^{\frac{1}{2}}=\\
&=-(\sigma+2b+1) - s (\sigma+1) + \\
&+ (1-s)\left[(\sigma-1)^2+ 4 \sigma r-4\sigma z+\rho^2 z^2
+ \rho^2 \left(y + \frac{b-1}{\sigma}x\right)^2\right]^{\frac{1}{2}}.
 \end{split}
\end{equation}

We can use the inequality
\begin{equation}\label{ineq:KL}
 \sqrt{k+l} \le \sqrt{k} + \frac{l}{2 \sqrt{k}}
\end{equation}
for  $k+l > 0, k > 0$.
Consider
$$k+l=(\sigma-1)^2+ 4 \sigma r-4\sigma z+\rho^2 z^2
+ \rho^2 \left(y + \frac{b-1}{\sigma}x\right)^2$$
and
 $$k=(\sigma-1)^2+4\sigma r.$$
According to \eqref{transformation} the following
equation
 \begin{equation*}
 \begin{split}
&k+l=(\sigma-1)^2+ 4 \sigma r-4\sigma z+\rho^2 z^2 + \rho^2
\left(y + \frac{b-1}{\sigma}x\right)^2=\\
&=(\sigma-2b+1)^2 + \left(\frac{2\sigma}{\rho} -\rho z \right)^2
+\rho^2 \left(y + \frac{b-1}{\sigma}x \right)^2 \ge 0
 \end{split}
 \end{equation*}
is valid.
By the condition of Theorem \ref{Main} we have
$$r\sigma+(\sigma-b)(b-1)>0 \Leftrightarrow r\sigma > -\sigma b +\sigma +b^2-b.$$
Then
$$k=(\sigma-1)^2+4\sigma r > (\sigma-1)^2 - 4\sigma b + 4\sigma +4b^2 - 4b = (\sigma - 2b + 1)^2 \ge 0.$$

Applying the inequality \eqref{ineq:KL} with $k=(\sigma-1)^2 + 4 \sigma r$ and
 $l=-4\sigma z+\rho^2 z^2 + \rho^2 \left(y + \frac{b-1}{\sigma}x\right)^2$
to the relation \eqref{transformation}, we obtain
\begin{equation}\label{lamda1lambda2slambda3}
 \begin{split}
&2 (\lambda_1+\lambda_2+s\lambda_3) \le -(\sigma + 2b + 1)-s(\sigma+1)
+(1-s)[(\sigma-1)^2+4\sigma r]^{\frac{1}{2}}+\\
&+\frac{2(1-s)}{[(\sigma-1)^2+4\sigma r]^{\frac{1}{2}}}
\left[-\sigma z+\frac{\rho^2 z^2}{4}+\frac{\rho^2}{4}
\left(y+\frac{b-1}{\sigma}x\right)^2\right].
 \end{split}
\end{equation}

Introduce the function
\begin{equation}\label{theta}
 \theta (x, y, z) = \frac{(1-s) V(x, y, z)}{[(\sigma-1)^2+4\sigma r]^{\frac{1}{2}}},
\end{equation}
where
\begin{equation}\label{V}
 V(x, y, z) = \gamma_4 x^2+(-\sigma\gamma_1+\gamma_3)y^2+\gamma_3z^2
+\frac{1}{4\sigma}\gamma_1 x^4-\gamma_1 x^2 z - \gamma_1 \gamma_2 xy -\frac{\sigma}{b}z.
\end{equation}
We choose the running parameters $\gamma_1, \gamma_2, \gamma_3, \gamma_4$
in such a way that
\begin{equation}\label{R}
 R := -\sigma z + \frac{\rho^2 z^2}{4} + \frac{\rho^2}{4}
\left(y + \frac{b-1}{\sigma}x\right)^2 + \dot V \le 0 \quad \forall x, y, z,
\end{equation}
where $\dot{V}(x)$ is a derivative with respect  to system \eqref{sys:Lorenz-classic}.
Then, using \eqref{lamda1lambda2slambda3}, we obtain
\begin{equation}\label{parametrs}
 2 (\lambda_1+\lambda_2+s\lambda_3) + 2\dot \theta \le
-(\sigma + 2b + 1)-s(\sigma+1)+(1-s)[(\sigma-1)^2+4\sigma r]^{\frac{1}{2}},
\end{equation}
i.e. we estimated $2 (\lambda_1+\lambda_2+s\lambda_3) + 2\dot \theta$
by the relation, which depends on the parameters of system \eqref{sys:Lorenz-classic}
and is independent of $x, y, z$.

\bigskip
\textbf{IV.}
We perform the analysis of $R$, choosing the running parameters
$\gamma_1, \gamma_2, \gamma_3, \gamma_4$ in such a way that \eqref{R} is valid.
Taking into account \eqref{V}, we obtain
$$R = -\gamma_1 x^4+(\gamma_1 \gamma_2+2\sigma\gamma_1+b\gamma_1)x^2 z
+\left(\frac{\rho^2(b-1)^2}{4\sigma^2}-r\gamma_1 \gamma_2-2\sigma\gamma_4\right)x^2+$$
$$+\left(\frac{\rho^2(b-1)}{2\sigma}+\gamma_1 \gamma_2-2r\sigma\gamma_1
+2r\gamma_3+2\sigma\gamma_4+\sigma\gamma_1 \gamma_2-\frac{\sigma}{b}\right)x y+$$
$$+\left(2\sigma\gamma_1-2\gamma_3-\sigma\gamma_1 \gamma_2
+\frac{\rho^2}{4}\right)y^2+\left(\frac{\rho^2}{4}-2b\gamma_3\right)z^2=$$
$$=A_1 x^4 + A_2x^2 z + A_3 z^2 + B_1 x^2 + B_2 x y + B_3 y^2 ,$$
where
$$A_1 = -\gamma_1, A_2 = \gamma_1(\gamma_2+2\sigma+b), A_3 = \frac{\rho^2}{4}-2b\gamma_3$$
$$B_1 = \frac{\rho^2(b-1)^2}{4\sigma^2}-r\gamma_1 \gamma_2-2\sigma\gamma_4
= C_1 - 2\sigma\gamma_4$$
$$B_2 = \frac{\rho^2(b-1)}{2\sigma}+\gamma_1 \gamma_2-2r\sigma\gamma_1
+2r\gamma_3+\sigma\gamma_1 \gamma_2-\frac{\sigma}{b}+2\sigma\gamma_4=C_2+2\sigma\gamma_4$$
$$B_3 = 2\sigma\gamma_1-2\gamma_3-\sigma \gamma_1 \gamma_2+\frac{\rho^2}{4}.$$

\bigskip
\textbf{IV. a.}
Under the conditions $A_1 \ne 0$ and $B_3 \ne 0$ we transform
the relation for $R$:
\begin{equation}\label{RCombined}
 R = A_1\left(x^2+\frac{A_2}{2A_1}z\right)^2
+B_3\left(y+\frac{B_2}{2B_3}x\right)^2+\frac{4B_1 B_3 - B_2^2}{4B_3}x^2
+\frac{4A_1 A_3 - A_2^2}{4A_1}z^2
\end{equation}

Then
\begin{equation} \label{RConditions}
 \left.
 \begin{aligned}
 &A_1 < 0\\
 &B_3 < 0\\
 &4A_1A_3-A_2^2 \ge 0\\
 &4B_1 B_3-B_2^2 \ge 0
 \end{aligned}
 \right\}
 \Rightarrow R \le 0 \quad \forall x, y, z
\end{equation}
Now we analyze the signs of four addends in \eqref{RCombined}

(1) \begin{equation}\label{ineq:A1final}
 A_1 = -\gamma_1 < 0 \Leftrightarrow \gamma_1>0
 \end{equation}

(2)
\begin{equation}\label{ineq:B3final}
 B_3 < 0 \Leftrightarrow
 2\gamma_3 > 2\sigma\gamma_1-\sigma\gamma_1\gamma_2+\frac{\rho^2}{4}
\end{equation}

(3) \begin{equation}\label{ineq:A1A3A2}
 0 \le 4A_1 A_3-A_2^2=-4\gamma_1\left(\frac{\rho^2}{4}-2b\gamma_3\right)
-{\gamma_1}^2(\gamma_2+2\sigma+b)^2.
\end{equation}
Since according to \eqref{ineq:A1final}, $\gamma_1 >0$,
from \eqref{ineq:A1A3A2} it follows that
$$-\gamma_1(\gamma_2+2\sigma+b)^2 -\rho^2+8b\gamma_3 \ge 0$$
$$\Leftrightarrow$$
\begin{equation}\label{ineq:A1A3A2final}
 2\gamma_3 \ge \frac{\rho^2}{4b} + \frac{\gamma_1}{4b}(\gamma_2 + 2 \sigma +b)^2.
\end{equation}

(4) $$4B_1 B_3-B_2^2=4B_3 (C_1-2\sigma\gamma_4)-(C_2+2\sigma\gamma_4)^2=$$
$$=4B_3 C_1-8\sigma B_3\gamma_4-C_2^2-4\sigma C_2\gamma_4-4\sigma^2\gamma_4^2=$$
$$=-4\sigma^2\gamma_4^2+4\gamma_4\sigma(-2B_3-C_2)+4B_3 C_1-C_2 ^2.$$
The relation $4B_1 B_3-B_2^2$ is a quadratic polynomial in $\gamma_4$
with a negative coefficient of $\gamma_4^2$. Therefore for the inequality
$4B_1 B_3-B_2^2 \ge 0$ to be satisfied, it is necessary that the corresponding
quadratic equation has a real root, i.e. a positive discriminant:
$$D_{\gamma_4} = 16\sigma^2(4B_3^2+C_2^2+4B_3C_2)+16\sigma^2(4B_3C_1-C_2^2)=$$
$$=16 C_2^2\sigma^2+64\sigma^2B_3^2+64\sigma^2B_3 C_2+64\sigma^2B_3C_1
-16\sigma^2C_2^2=64\sigma^2B_3(B_3+C_2+C_1).$$
Since according to \eqref{RConditions} the condition $B_3 < 0$ is satisfied,
we obtain $D_{\gamma_4} \ge 0 \Leftrightarrow B_3+C_1+C_2 \le 0$. We have
$$B_3+ÔC_1+C_2=2\sigma\gamma_1-2\gamma_3-\sigma\gamma_1\gamma_2
+\frac{\rho^2}{4}+\frac{\rho^2(b-1)^2}{4\sigma^2}+$$
$$-r\gamma_1\gamma_2+\frac{\rho^2(b-1)}{2\sigma}+\gamma_1\gamma_2
-2r\sigma\gamma_1+2r\gamma_3+\sigma\gamma_1\gamma_2-\frac{\sigma}{b}=$$
$$=\gamma_3(2r-2)-\gamma_1(-2\sigma+2r\sigma)
-\gamma_1\gamma_2(\sigma+r-\sigma-1)+\frac{\rho^2}{4}+\frac{\rho^2(b-1)^2}{4\sigma^2}+$$
$$+\frac{\rho^2(b-1)}{2\sigma}-\frac{\sigma}{b}
=2\gamma_3(r-1)-2\sigma\gamma_1(r-1)-\gamma_1\gamma_2(r-1)+\frac{\rho^2}{4}+$$
$$+\frac{\rho^2(b-1)^2}{4\sigma^2}+\frac{\rho^2(b-1)}{2\sigma}-\frac{\sigma}{b}.$$
Thus, for
\begin{equation}\label{ineq:B1B3B2final}
 \begin{split}
 &2\gamma_3(r-1) \le 2\sigma\gamma_1(r-1) + \gamma_1\gamma_2(r-1)
-\frac{\rho^2}{4}-\\
 &-\frac{\rho^2(b-1)^2}{4\sigma^2}-\frac{\rho^2(b-1)}{2\sigma}
+\frac{\sigma}{b}
 \end{split}
\end{equation}
there exists $\gamma_4$ such that $4B_1 B_3-B_2^2 \ge 0$.

Note that a part of the right-hand side of the inequality in \eqref{ineq:B1B3B2final}
can be transform in the following way
\begin{equation}\label{transform-to-complete-square}
 \begin{split}
 &\frac{\rho^2}{4}+\frac{\rho^2(b-1)^2}{4\sigma^2}
+\frac{\rho^2(b-1)}{2\sigma}=\\
 &=\frac{\rho^2}{\sigma^2}\left(\frac{\sigma^2}{4}
+\frac{(b-1)^2}{4}+\frac{\sigma(b-1)}{2}\right)=\frac{\rho^2(\sigma+b-1)^2}{4\sigma^2}.
 \end{split}
\end{equation}

Taking into account \eqref{ineq:A1final}, \eqref{ineq:B3final}, \eqref{ineq:A1A3A2final},
\eqref{ineq:B1B3B2final}, \eqref{transform-to-complete-square},
conditions \eqref{RConditions} become
\begin{equation}\label{RConditionsA1B3}
 \left.
 \begin{aligned}
 &\gamma_1 > 0\\
 &2\gamma_3 > 2\sigma\gamma_1-\sigma \gamma_1 \gamma_2+\frac{\rho^2}{4}\\
 &2\gamma_3 \ge \frac{\rho^2}{4b} + \frac{\gamma_1}{4b}(\gamma_2
+ 2 \sigma +b)^2\\
 &2(r-1)\gamma_3 \le 2\sigma\gamma_1(r-1)
+ \gamma_1\gamma_2(r-1)-\frac{\rho^2(b+\sigma-1)^2}{4\sigma^2}+\frac{\sigma}{b}\\
 \end{aligned}
 \right \}
 \Rightarrow R \le 0 \quad \forall x, y, z.
\end{equation}

For obtaining \eqref{RCombined} it is assumed that $A_1 \ne 0, B_3 \ne 0$.
Let us analyze the cases $A_1 = 0$ and $B_3=0$.

\bigskip
\textbf{IV. b.}
Consider $B_3 =0$. In this case we have
$$R = A_1 x^4 + A_2x^2 z + A_3 z^2 + B_1 x^2 + B_2 x y,$$ i.e. it is impossible
to choose $\gamma_1, \gamma_2, \gamma_3, \gamma_4$
in such a way that $R \le 0$ is valid $\forall x, y, z$.

\bigskip
\textbf{IV. c.}
Consider the case $\gamma_1=A_1=0$ and $B_3 \ne 0$.
Then
$$ A_2 = \gamma_1(\gamma_2+2\sigma+b) =0$$
and
$$R = A_3 z^2 + B_1 x^2 + B_2 x y + B_3 y^2 = $$
$$=A_3 z^2+B_3\left(y+\frac{B_2}{2B_3}x\right)^2+\frac{4B_1 B_3 - B_2^2}{4B_3}x^2.$$
In this case
\begin{equation}
 \left.
 \begin{aligned}
 &A_3 \le 0\\
 &B_3 < 0\\
 &4B_1 B_3-B_2^2 \ge 0
 \end{aligned}
 \right\}
 \Rightarrow R \le 0 \quad \forall x, y, z.
\end{equation}
The second and third conditions are similar to the second and fourth conditions
in \eqref{RConditions}. Consequently it remains to consider
$$0 \ge A_3= \frac{\rho^2}{4}-2b\gamma_3
\Leftrightarrow
2\gamma_3 \ge \frac{\rho^2}{4b}.$$
The latter inequality, obtained under the assumption $A_1=0$, coincides
with condition \eqref{ineq:B3final} if in \eqref{ineq:B3final}
we consider $\gamma_1 = 0$. Thus, the conditions
on the running parameters $\gamma_1, \gamma_2, \gamma_3, \gamma_4$,
obtained under the condition $A_1 = 0, B_3 \ne 0$, can be joined
with those, obtained under the condition $A_1 \ne 0, B_3 \ne 0$, i.e.
\begin{equation}\label{RConditionsFinal}
 \left.
 \begin{aligned}
 &\gamma_1 \ge 0\\
 &2\gamma_3 > 2\sigma\gamma_1-\sigma \gamma_1 \gamma_2+\frac{\rho^2}{4}\\
 &2\gamma_3 \ge \frac{\rho^2}{4b} + \frac{\gamma_1}{4b}(\gamma_2
+ 2 \sigma +b)^2\\
 &2(r-1)\gamma_3 \le 2\sigma\gamma_1(r-1)
+ \gamma_1\gamma_2(r-1)-\frac{\rho^2(b+\sigma-1)^2}{4\sigma^2}+\frac{\sigma}{b}\\
 \end{aligned}
 \right\}
 \Rightarrow
 R \le 0 \quad \forall x, y, z.
\end{equation}

\bigskip
\textbf{V.}
For  $r-1>0$ conditions \eqref{RConditionsFinal} can be transformed
in the following way:
\begin{equation}
 \left\{
 \begin{aligned}
 &\gamma_1 \ge 0,\\
 &2\gamma_3 > \frac{\rho^2}{4} + 2\sigma \gamma_1
- \sigma \gamma_1 \gamma_2,\\
 &2\gamma_3 \ge \frac{\rho^2}{4b} + \frac{\gamma_1}{4b}
(\gamma_2 + 2 \sigma +b)^2,\\
 &2\gamma_3 \le 2\sigma \gamma_1 + \gamma_1 \gamma_2
- \frac{\rho^2(b+\sigma -1)^2}{4\sigma^2 (r-1)} + \frac{\sigma}{b(r-1)}.
 \end{aligned}
 \right.
\end{equation}

For the existence of $\gamma_3$ it is necessary and sufficient that
\begin{equation}\label{gamma3existscond}
 \left\{
 \begin{gathered}
 \frac{\rho^2}{4} + \gamma_1(2 \sigma - \sigma \gamma_2)
< \gamma_1(2 \sigma + \gamma_2) - \frac{\rho^2(b + \sigma - 1)^2}
{4 \sigma^2(r-1)} + \frac{\sigma}{b(r-1)}, \hfill \\
 \frac{\rho^2}{4b} + \frac{\gamma_1}{4b}(2\sigma + b
+ \gamma_2)^2 \le \gamma_1(2 \sigma + \gamma_2)
- \frac{\rho^2(b + \sigma - 1)^2}{4 \sigma^2(r-1)} + \frac{\sigma}{b(r-1)}. \hfill \\
 \end{gathered}
 \right.
\end{equation}

We now analyze the obtained inequalities.

\bigskip
\textbf{V. a.}
Consider the first inequality from system \eqref{gamma3existscond}
for different $\gamma_2$.

Let be $\gamma_2=0$. Then the inequality is equivalent to
$$\frac{\rho^2}{4} + \frac{\rho^2(b+\sigma-1)^2}{4\sigma^2(r-1)}
- \frac{\sigma}{b(r-1)} < 0.$$

If $\gamma_2 > 0$, then we obtain the condition
$$\gamma_1 > \frac{1}{\gamma_2(\sigma+1)}\left(\frac{\rho^2}{4}
+ \frac{\rho^2(b+\sigma-1)^2}{4\sigma^2(r-1)} - \frac{\sigma}{b(r-1)}\right),$$
and if $\gamma_2 < 0$, then
$$\gamma_1 < \frac{1}{\gamma_2(\sigma+1)}\left(\frac{\rho^2}{4}
+ \frac{\rho^2(b+\sigma-1)^2}{4\sigma^2(r-1)} - \frac{\sigma}{b(r-1)}\right).$$
In the latter case it is required that
\begin{equation}\label{ineq:gamma1ExistGamma2Neg}
 \frac{\rho^2}{4} + \frac{\rho^2(b+\sigma-1)^2}{4\sigma^2(r-1)}
- \frac{\sigma}{b(r-1)} < 0
\end{equation}
since according to \eqref{RConditionsFinal} we have $\gamma_1 \ge 0$.
Thus, since the cases of $\gamma_2=0$ and $\gamma_2<0$
impose the same condition on the parameters of the system, they can be joined.

\bigskip
\textbf{V. b.}
From the second inequality of system \eqref{gamma3existscond}
we obtain condition for $\gamma_1$:
$$\frac{(2\sigma - b + \gamma_2)^2}{4b} \gamma_1 \le
-\left(\frac{\rho^2}{4b} + \frac{\rho^2(b+\sigma-1)^2}{4\sigma^2(r-1)}
- \frac{\sigma}{b(r-1)}\right).$$

If $\gamma_2 = b-2\sigma$, then it is required that
$$\frac{\rho^2}{4b} + \frac{\rho^2(b+\sigma-1)^2}{4\sigma^2(r-1)}
- \frac{\sigma}{b(r-1)} \le 0.$$
If $\gamma_2 \ne b-2\sigma$, then we obtain the following condition for $\gamma_1$
$$\gamma_1 \le -\frac{4b}{(2\sigma - b + \gamma_2)^2}\left(\frac{\rho^2}{4b}
+ \frac{\rho^2(b+\sigma-1)^2}{4\sigma^2(r-1)} - \frac{\sigma}{b(r-1)}\right).$$
According to \eqref{RConditionsFinal}, it is valid $\gamma_1\ge0$.
Therefore it is required the following
\begin{equation}\label{ineq:gamma1Exist}
 \frac{\rho^2}{4b} + \frac{\rho^2(b+\sigma-1)^2}{4\sigma^2(r-1)}
- \frac{\sigma}{b(r-1)} \le 0.
\end{equation}
Thus, inequality \eqref{ineq:gamma1Exist}
holds true also for  $\gamma_2 = b-2\sigma$ and for  $\gamma_2 \ne b-2\sigma$.

Condition \eqref{ineq:gamma1Exist} must be satisfied
for any sign of $\gamma_2$. If, in addition, condition
\eqref{ineq:gamma1ExistGamma2Neg} is also satisfied,
then in function \eqref{V}
we can take $\gamma_2 \le 0$ and in this case there exist
$\gamma_1, \gamma_3, \gamma_4$ such that the relation
$$R = -\sigma z + \frac{\rho^2 z^2}{4} + \frac{\rho^2}{4}
\left(y + \frac{b-1}{\sigma}x\right)^2 + \dot V \le 0 \quad \forall x, y, z $$
is valid.
If inequality \eqref{ineq:gamma1ExistGamma2Neg} is not satisfied,
then in function \eqref{V} it must be considered $\gamma_2>0$ and,
except for condition \eqref{ineq:gamma1Exist}, it is necessary to find
conditions for the existence of $\gamma_1$.

In this case we obtain the following family:

\begin{equation}\label{d1Interim}
 \left[
 \begin{gathered}
 \left\{
 \begin{gathered}
 \frac{\rho^2}{4b} + \frac{\rho^2(b+\sigma-1)^2}{4\sigma^2(r-1)}
- \frac{\sigma}{b(r-1)} \le 0, \hfill \\
 \frac{\rho^2}{4} + \frac{\rho^2(b+\sigma-1)^2}{4\sigma^2(r-1)}
- \frac{\sigma}{b(r-1)} < 0 \hfill \\
 \end{gathered}
 \right. \hfill\\
 \left\{
 \begin{gathered}
 \frac{\rho^2}{4b} + \frac{\rho^2(b+\sigma-1)^2}{4\sigma^2(r-1)}
- \frac{\sigma}{b(r-1)} \le 0 , \hfill \\
 \gamma_2 > 0, \gamma_2 \ne b-2\sigma, \hfill \\
 \frac{1}{\gamma_2(\sigma+1)}\left(\frac{\rho^2}{4}
+ \frac{\rho^2(b+\sigma-1)^2}{4\sigma^2(r-1)} - \frac{\sigma}{b(r-1)}\right) < \hfill \\
 <-\frac{4b}{(2\sigma - b + \gamma_2)^2}\left(\frac{\rho^2}{4b}
+ \frac{\rho^2(b+\sigma-1)^2}{4\sigma^2(r-1)} - \frac{\sigma}{b(r-1)}\right).
 \\
 \end{gathered}
 \right. \hfill\\
 \end{gathered}
 \right.
\end{equation}

\bigskip
\textbf{V. c.}

Consider the first condition of family \eqref{d1Interim}:
$$\left\{ \begin{gathered}
 \frac{\rho^2}{4b} + \frac{\rho^2(b+\sigma-1)^2}{4\sigma^2(r-1)}
- \frac{\sigma}{b(r-1)} \le 0, \hfill \\
 \frac{\rho^2}{4} + \frac{\rho^2(b+\sigma-1)^2}{4\sigma^2(r-1)}
- \frac{\sigma}{b(r-1)} < 0.\hfill \\
\end{gathered} \right.$$

Substituting the value $\rho$ from \eqref{rho}, we obtain
$$\left\{ \begin{gathered}
 \frac{\sigma^2}{4b((r-1)\sigma-b^2+b+\sigma b)}
+ \frac{(b+\sigma-1)^2}{4(r-1)((r-1)\sigma-b^2+b+\sigma b)}
- \frac{\sigma}{b(r-1)} \le 0, \hfill \\
 \frac{\sigma^2}{4((r-1)\sigma-b^2+b+\sigma b)}
+ \frac{(b+\sigma-1)^2}{4(r-1)((r-1)\sigma-b^2+b+\sigma b)}
- \frac{\sigma}{b(r-1)} < 0.
\end{gathered} \right.$$
$$\Leftrightarrow$$
$$\left\{ \begin{gathered}
 r-1 \ge \frac{b(b+\sigma-1)^2 - 4\sigma(b+\sigma b -b^2)}{3\sigma^2}, \hfill \\
 \sigma^2 (r-1)(b-4) < 4\sigma (\sigma b + b -b^2) - b(b+\sigma - 1)^2.
\end{gathered} \right.$$

\bigskip
\textbf{V. d.}
Consider the second condition of family \eqref{d1Interim}:
\begin{equation}\label{d1Interim2}
 \left\{
 \begin{gathered}
 \frac{\rho^2}{4b} + \frac{\rho^2(b+\sigma-1)^2}{4\sigma^2(r-1)}
- \frac{\sigma}{b(r-1)} \le 0 , \hfill \\
 \gamma_2 > 0, \gamma_2 \ne b-2\sigma, \hfill \\
 \frac{1}{\gamma_2(\sigma+1)}\left(\frac{\rho^2}{4}
+ \frac{\rho^2(b+\sigma-1)^2}{4\sigma^2(r-1)} - \frac{\sigma}{b(r-1)}\right) < \hfill \\
 <-\frac{4b}{(2\sigma - b + \gamma_2)^2}\left(\frac{\rho^2}{4b}
+ \frac{\rho^2(b+\sigma-1)^2}{4\sigma^2(r-1)} - \frac{\sigma}{b(r-1)}\right).
 \end{gathered}
 \right.
\end{equation}

From Item \textbf{V. c.} it is known that the first inequality of system
is equivalent to the relation
$$
  r-1 \ge \frac{b(b+\sigma-1)^2 - 4\sigma(b+\sigma b -b^2)}{3\sigma^2}.
$$
Next we transform the inequality
$$\frac{1}{\gamma_2(\sigma+1)}\left(\frac{\rho^2}{4}
+ \frac{\rho^2(b+\sigma-1)^2}{4\sigma^2(r-1)} - \frac{\sigma}{b(r-1)}\right) < $$
$$<-\frac{4b}{(2\sigma - b + \gamma_2)^2}\left(\frac{\rho^2}{4b}
+ \frac{\rho^2(b+\sigma-1)^2}{4\sigma^2(r-1)} - \frac{\sigma}{b(r-1)}\right)$$
$$\Leftrightarrow$$
$$(2\sigma - b + \gamma_2)^2 \left(\frac{\rho^2}{4}
+ \frac{\rho^2(b+\sigma-1)^2}{4\sigma^2(r-1)} - \frac{\sigma}{b(r-1)}\right) + $$
$$+ 4b\gamma_2(\sigma+1)\left(\frac{\rho^2}{4b}
+ \frac{\rho^2(b+\sigma-1)^2}{4\sigma^2(r-1)} - \frac{\sigma}{b(r-1)}\right) <0$$
$$\Leftrightarrow$$
\begin{equation}\label{gamma2Inequality}
 \begin{split}
 (2\sigma - b + \gamma_2)^2 \left[b(b+\sigma-1)^2
- 4\sigma(\sigma b + b-b^2) + \sigma^2(r-1)(b-4)\right] + \\
 + 4b\gamma_2(\sigma+1)\left(b(b+\sigma-1)^2
- 4\sigma(\sigma b + b-b^2) - 3 \sigma^2(r-1)\right) <0.
 \end{split}
\end{equation}

The left-hand side of the inequality
is a quadratic polynomial in $\gamma_2$.
If coefficient of $\gamma_2^2$ is negative,
i.e.
$$b(b+\sigma-1)^2 - 4\sigma(\sigma b + b-b^2) + \sigma^2(r-1)(b-4) <0,$$
then there exists $\gamma_2>0$ satisfying inequality \eqref{gamma2Inequality}.

If the coefficient is equal to $0$, i.e.
$$b(b+\sigma-1)^2 - 4\sigma(\sigma b + b-b^2) + \sigma^2(r-1)(b-4) =0,$$
then the relation
$$4b\gamma_2(\sigma+1)\left[b(b+\sigma-1)^2 - 4\sigma(\sigma b + b-b^2)
- 3 \sigma^2(r-1)\right] <0$$
must be valid.
By \eqref{d1Interim2}, $\gamma_2>0$ is satisfied
and
$b(b+\sigma-1)^2 - 4\sigma(b+\sigma b -b^2)-3\sigma^2(r-1) \le 0,$
Then the inequality
$$4b\gamma_2(\sigma+1)\left[b(b+\sigma-1)^2 - 4\sigma(\sigma b + b-b^2)
- 3 \sigma^2(r-1)\right] <0$$
is satisfied.

If  the coefficient is positive, i.e.
$$b(b+\sigma-1)^2 - 4\sigma(\sigma b + b-b^2) + \sigma^2(r-1)(b-4) >0,$$
then for the existence of $\gamma_2>0$, satisfying \eqref{gamma2Inequality},
it is necessary and sufficient that there exist two different real roots
of the equation
\begin{equation}
 \begin{split}
 (2\sigma - b + \gamma_2)^2 \left[b(b+\sigma-1)^2
- 4\sigma(\sigma b + b-b^2) + \sigma^2(r-1)(b-4)\right] + \\
 + 4b\gamma_2(\sigma+1)\left(b(b+\sigma-1)^2
- 4\sigma(\sigma b + b-b^2) - 3 \sigma^2(r-1)\right) =0
 \end{split}
\end{equation}
and the largest root must be positive.
This condition corresponds to the second condition, stated in the theorem.

\bigskip
\textbf{V. e.}
Thus, there are obtained conditions \eqref{cond:mainTheoremR}-\eqref{cond:mainTheorem2}
for  which
$$R = -\sigma z + \frac{\rho^2 z^2}{4} + \frac{\rho^2}{4}
\left(y + \frac{b-1}{\sigma}x\right)^2 + \dot V \le 0 \quad \forall x, y, z .$$

Consequently, it is valid \eqref{parametrs}:
$$2 (\lambda_1+\lambda_2+s\lambda_3) + 2\dot \theta \le -(\sigma + 2b + 1)
-s(\sigma+1)+(1-s)[(\sigma-1)^2+4\sigma r]^{\frac{1}{2}}.$$

\bigskip
\textbf{VI.}
To apply Theorem \ref{theorem:th2}, we consider inequality \eqref{parametrs}
for  $s=0$:
$$2 (\lambda_1+\lambda_2 + \dot{\theta}) \le -(\sigma + 2b + 1)
+[(\sigma-1)^2+4\sigma r]^{\frac{1}{2}} .$$
Then we find the conditions for the parameters of system $\sigma, b, r$
for  which $2 (\lambda_1+\lambda_2 + \dot{\theta}) < 0$. We have
\begin{equation}\label{ineq:lamda1lambda2}
 -(\sigma + 2b + 1)+[(\sigma-1)^2+4\sigma r]^{\frac{1}{2}} < 0 \Leftrightarrow
[(\sigma-1)^2+4\sigma r]^{\frac{1}{2}} < (\sigma + 2b + 1).
\end{equation}
The parameters of system $\sigma, b$ are positive numbers, i.e. $\sigma + 2b + 1 > 0$.
Therefore after squaring two sides of inequality \eqref{ineq:lamda1lambda2}
we obtain
$$\sigma^2 - 2\sigma + 1 + 4 \sigma r < \sigma^2 + 1 + 4 b^2 + 2 \sigma + 4 \sigma b + 4 b$$
$$\Leftrightarrow$$
$$4\sigma r < 4b^2 + 4\sigma + 4\sigma b + 4b$$
$$\Leftrightarrow$$
$$\sigma r < (b + 1)(b + \sigma)$$
Thus, taking into account inequality \eqref{ineq:sqrt},
we obtain condition \eqref{cond:mainTheoremR1} for  which
$$2 (\lambda_1+\lambda_2 + \dot{\theta}) < 0,$$
and Theorem \ref{theorem:th2} can be applied,
i.e. any bounded on $[0; +\infty)$ solution of system \eqref{sys:Lorenz-classic}
tends to a certain equilibrium as $t \to +\infty$.

\bigskip
\textbf{VII.}
To apply Theorem \ref{theorem:th1}, we consider $s \ne 0$ and \eqref{parametrs}:
$$2 (\lambda_1+\lambda_2+s\lambda_3 + \dot{\theta})
\le -(\sigma + 2b + 1)-s(\sigma+1)+(1-s)[(\sigma-1)^2+4\sigma r]^{\frac{1}{2}} =$$
$$= -s(\sigma+1 + [(\sigma-1)^2+4\sigma r]^{\frac{1}{2}}) -(\sigma + 2b + 1)
+ [(\sigma-1)^2+4\sigma r]^{\frac{1}{2}}.$$
We find $s$ such that the relation $2 (\lambda_1+\lambda_2+s\lambda_3
+ \dot{\theta}) < 0$ is valid:
$$-s(\sigma+1 + [(\sigma-1)^2+4\sigma r]^{\frac{1}{2}}) -(\sigma + 2b + 1)
 + [(\sigma-1)^2+4\sigma r]^{\frac{1}{2}} < 0$$
$$\Leftrightarrow$$
\begin{equation}\label{SInequality}
 -(\sigma + 2b + 1) + [(\sigma-1)^2+4\sigma r]^{\frac{1}{2}} <
s(\sigma+1 + [(\sigma-1)^2+4\sigma r]^{\frac{1}{2}}).
\end{equation}

The parameter $\sigma$ of system \eqref{sys:Lorenz-classic}
is a positive number.
Therefore the coefficient of $s$,
equal to $\sigma+1 + [(\sigma-1)^2+4\sigma r]^{\frac{1}{2}}$,
is greater than $0$.
In this case if inequality \eqref{SInequality} is divided
by $\sigma+1 + [(\sigma-1)^2+4\sigma r]^{\frac{1}{2}}$,
then the sign is not changed, i.e. we have
\begin{equation}\label{condS}
 s > s_0=\frac{-(\sigma + 2b + 1)
+[(\sigma-1)^2+4\sigma r]^{\frac{1}{2}}}{\sigma + 1
+ [(\sigma-1)^2+4\sigma r]^{\frac{1}{2}}}.
\end{equation}
Thus, in the case when \eqref{cond:mainTheoremR}-\eqref{cond:mainTheorem2}
is satisfied and $s$ satisfies \eqref{condS}
we have $2 (\lambda_1+\lambda_2+s\lambda_3 + \dot{\theta} ) < 0$.

Under the hypothesis of Theorem \ref{theorem:th1} we have $s \in [0,1)$.
However according to \eqref{condS} the relation $s > s_0 $ must be valid.
The case $s=0$ is already considered. Obviously
$\sigma+ 1 + [(\sigma-1)^2+4\sigma r]^{\frac{1}{2}} >0$.
Therefore for $s_0>0$ to be valid, it is required that
$$-(\sigma + 2b + 1)+[(\sigma-1)^2+4\sigma r]^{\frac{1}{2}} >0$$
$$\Leftrightarrow$$
\begin{equation}\label{ineq:lambda1lamba2slambda3}
 [(\sigma-1)^2+4\sigma r]^{\frac{1}{2}} > \sigma + 2b + 1.
\end{equation}
The relation $\sigma + 2b +1 \ge 0$ is always satisfied
since the parameters $\sigma, b$ of system \eqref{sys:Lorenz-classic}
are positive numbers.
Thus, after squaring inequality \eqref{ineq:lambda1lamba2slambda3}
we obtain
$$\sigma^2 - 2\sigma + 1 + 4\sigma r > \sigma^2 + 1 + 4b^2 + 2\sigma + 4\sigma b + 4b$$
$$\Leftrightarrow$$
\begin{equation}\label{ineq:sigmaRToCompare}
 \sigma r > (b + 1)(b + \sigma).
\end{equation}
Next we find a condition for  which $s_0 < 1$:
$$\frac{-(\sigma + 2b + 1)+[(\sigma-1)^2+4\sigma r]^{\frac{1}{2}}}{\sigma + 1
+ [(\sigma-1)^2+4\sigma r]^{\frac{1}{2}}} < 1$$
$$\Leftrightarrow$$
$$-(\sigma + 2b + 1)+[(\sigma-1)^2+4\sigma r]^{\frac{1}{2}} < \sigma + 1
+ [(\sigma-1)^2+4\sigma r]^{\frac{1}{2}}$$
$$\Leftrightarrow$$
$$0 < \sigma + 1 + b$$

The latter inequality is always valid. The parameters of system \eqref{sys:Lorenz-classic}
are positive numbers and therefore
$$(b+1)(b+\sigma)>(b-\sigma)(b-1)$$ is always satisfied,
i.e. from conditions \eqref{ineq:sqrt}, \eqref{ineq:sigmaRToCompare}
it remains only condition \eqref{ineq:sigmaRToCompare}.
Thus, it is obtained condition \eqref{cond:mainTheoremR2}
under which the relation $2 (\lambda_1+\lambda_2+s\lambda_3 + \dot{\theta} ) < 0$
is satisfied and Theorem \ref{theorem:th1} can be applied.
Consequently, $\dim_L K \le 2 + s$ for  $s$, satisfying \eqref{condS}
and \eqref{cond:mainTheoremR2}. Thus,
$$\dim_L K \le 3 - \frac{2 (\sigma + b + 1)}{\sigma + 1 + \sqrt{(\sigma-1)^2
+ 4 \sigma r}}.$$
\end{proof}

\begin{proof}(Lemma 1)
For positive parameters we have
$$(b+1)(b+\sigma)>(b-\sigma)(b-1).$$
Consequently if condition \eqref{cond:mainTheoremR2} is satisfied,
then the inequality
$$\sigma r > (b-1)(b-\sigma)$$ is satisfied too.
First we transform this inequality
$$r\sigma + (\sigma - b)(b - 1) > 0 \Leftrightarrow b^2 - b(\sigma+1)
+ 1 \sigma - r\sigma <0$$
$$ \Leftrightarrow$$
\begin{equation} \label{ineq:b12}
 b_1 < b < b_2,
\end{equation}
where $b_{1,2} = \frac{(\sigma+1) \mp \sqrt{(\sigma-1)^2 + 4 r\sigma}}{2}$.

When considered the linearized system
along the zero solution the corresponding matrix $J({\bf \rm x})=J$
is constant:
$$J = \begin{pmatrix}
-\sigma & \sigma & 0 \\
r & -1 & 0 \\
0 & 0 & -b
\end{pmatrix}.$$
Its eigenvalues are real and have the form:
\begin{align*}
 &\lambda_1 = -\frac{1}{2}(\sigma + 1 - \sqrt{(\sigma-1)^2 + 4r\sigma})\\
 &\lambda_2 = -b\\
 &\lambda_3 = -\frac{1}{2}(\sigma + 1 + \sqrt{(\sigma-1)^2 + 4r\sigma}).
\end{align*}
Then $\LE_i(0,0,0) = \lambda_i, i=1,2,3$.

 Condition \eqref{ineq:b12} is valid:
$$b_2 > b > b_1 \Leftrightarrow -\frac{1}{2}(\sigma+1 - \sqrt{(\sigma-1)^2
+ 4r\sigma}) > -b > -\frac{1}{2}(\sigma+1 + \sqrt{(\sigma-1)^2 + 4r\sigma}).$$
Thus, $\LE_1(0,0,0) > \LE_2(0,0,0) > \LE_3(0,0,0)$.

We shall show that if \eqref{cond:mainTheoremR2} is satisfied
and $j(t, (0,0,0))=2$, then
the definition of \eqref{formula:kaplan} satisfies all requirements of
the definition of local Lyapunov dimension.

(1) $$\LE_1(0,0,0) + \LE_2(0,0,0) = -\frac{1}{2}(\sigma + 1 - \sqrt{(\sigma-1)^2 + 4r\sigma}) - b =$$
$$=-\frac{1}{2}(\sigma+2b+1) + \frac{1}{2}\sqrt{(\sigma-1)^2 + 4r\sigma}.$$
The relation $\LE_1(0,0,0) + \LE_2(0,0,0) >0$ must be valid:
$$-\frac{1}{2}(\sigma+2b+1) + \frac{1}{2}\sqrt{(\sigma-1)^2 + 4r\sigma} > 0$$
$$\Leftrightarrow$$
$$\sqrt{(\sigma-1)^2 + 4r\sigma} > \sigma+2b+1$$
We obtained inequality \eqref{ineq:lambda1lamba2slambda3},
which, as is known, is equivalent to inequality \eqref{cond:mainTheoremR2}.

(2)
$$\LE_3(0,0,0) = -\frac{1}{2}(\sigma + 1 + \sqrt{(\sigma-1)^2 + 4r\sigma}) < -b.$$
Parameter $b$ of system \eqref{sys:Lorenz-classic} is positive,
consequently $\LE_3(0,0,0)<0$ is satisfied.

(3)
$$\frac{\LE_1(0,0,0) + \LE_2(0,0,0)}{|\LE_3(0,0,0)|}
= \frac{-\frac{1}{2}(\sigma + 1 - \sqrt{(\sigma-1)^2 + 4r\sigma}) - b}{|-\frac{1}{2}
(\sigma + 1 + \sqrt{(\sigma-1)^2 + 4r\sigma})|} = $$
$$=\frac{-(\sigma+2b+1)+\sqrt{(\sigma-1)^2+4r\sigma}}
{\sigma+1+\sqrt{(\sigma-1)^2+4r\sigma}}$$
The expression obtained coincides with $s_0$ of inequality \eqref{condS}.
The relation
$$\frac{\LE_1(0,0,0) + \LE_2(0,0,0)}{|\LE_3(0,0,0)|}=s_0 <1$$ must be valid.
According to the proof of Theorem \ref{Main}, we have $s_0 <1 \Leftrightarrow \sigma+1+b>0$.
This is true since the parameters of system $\sigma, b$ are positive numbers.

Thus, we have
$$\dim_L(0,0,0)=2 - \frac{\sigma +2b + 1 - \sqrt{(\sigma-1)^2
+4r\sigma}}{\sigma+1 + \sqrt{(\sigma-1)^2 + 4r\sigma}}
= 3 - \frac{2(\sigma +b+1)}{\sigma+1 + \sqrt{(\sigma-1 )^2 + 4r\sigma}}.$$
\end{proof}

\begin{proof} (Lemma 2)

 \textbf{I.}
 We consider sufficient condition for the theorem to be valid:
 \begin{equation}\label{cond:sufficientLemma2}
 \left\{
 \begin{gathered}
 r - 1 > \frac{b(b + \sigma - 1)^2 - 4\sigma(b + \sigma b
- b^2)}{3\sigma^2} \hfill \\
 \sigma^2 (r - 1)(b - 4) < 4\sigma (\sigma b + b - b^2)
- b(b + \sigma - 1)^2 \hfill \\
 \end{gathered}
 \right. \hfill\\
 \end{equation}
 and show that the condition is valid for sufficiently great $\sigma$.
 The domain in which condition \eqref{cond:sufficientLemma2}
is not satisfied is bounded with respect to $\sigma$.
 Since it is considered $b < 4, \sigma > 0$, system \eqref{cond:sufficientLemma2}
 is equivalent to the following system
 \begin{equation*}
 \left\{
 \begin{gathered}
 r - 1 > \frac{b(b + \sigma - 1)^2
- 4\sigma(b + \sigma b - b^2)}{3\sigma^2} \hfill \\
 r - 1 > \frac{b(b + \sigma - 1)^2
- 4\sigma(b + \sigma b - b^2)}{\sigma^2 (4-b)} \hfill \\
 \end{gathered}
 \right. \hfill\\
 \end{equation*}
 \begin{align*}\Leftrightarrow\end{align*}
 \begin{equation}\label{cond:afterDivision}
 \left\{
 \begin{gathered}
 r - 1 > \frac{b\left(\frac{b}{\sigma} + 1
- \frac{1}{\sigma}\right)^2 - 4\left(\frac{b}{\sigma}
+ b - \frac{b^2}{\sigma}\right)}{3} \hfill \\
 r - 1 > \frac{b\left(\frac{b}{\sigma} + 1
- \frac{1}{\sigma}\right)^2 - 4\left(\frac{b}{\sigma}
+ b - \frac{b^2}{\sigma}\right)}{4-b} \hfill \\
 \end{gathered}
 \right. \hfill\\
 \end{equation}
 Let be
 \begin{align}\label{functionF}
 f(b,\sigma) = \left(\frac{b}{\sigma} + 1
- \frac{1}{\sigma}\right)^2 - 4\left(\frac{1}{\sigma} + 1
- \frac{b}{\sigma}\right) = \hfill \\
 =\frac{1}{\sigma^2} b^2
+\frac{-2+6\sigma}{\sigma^2}b+\frac{-6\sigma+1-3\sigma^2}{\sigma^2}
 \end{align}
 Then system \eqref{cond:afterDivision} can be represented as
 \begin{equation}\label{cond:afterReplacement}
 \left\{
 \begin{gathered}
 r - 1 > \frac{b}{3} f(b,\sigma) \hfill \\
 r - 1 > \frac{b}{4-b} f(b,\sigma) \hfill \\
 \end{gathered}
 \right. \hfill\\
 \end{equation}
The function $f(b,\sigma)$ is quadratic polynomial in $b$
with a positive coefficient of $b^2$.
The abscissa of parabola peak is equal to $b_0 = 1 - 3 \sigma$.
If we consider $\sigma > \frac{1}{3}$,
then $b_0<0$ and, therefore, on the interval $b \in (0, 4)$
the maximum of function $f(b,\sigma)$
is attained for  $b=4$ and is equal to
\begin{equation*}
 f(4,\sigma) = -3+\frac{18}{\sigma}+ \frac{9}{\sigma^2}.
\end{equation*}
For  $\sigma > 3 + \sqrt{12}$ the value $f(4,\sigma) < 0$.

Since $b \in (0, 4)$ and $r>1$, for sufficiently great $\sigma$ we have
\begin{equation}
\left\{
\begin{gathered}
r - 1 > 0 > \frac{b}{3} f(b,\sigma) \hfill \\
r - 1 > 0 > \frac{b}{4-b} f(b,\sigma), \hfill \\
\end{gathered}
\right. \hfill\\
\end{equation}
i.e. condition \eqref{cond:sufficientLemma2} is satisfied.

\textbf{II.}
Let us show that for  sufficiently great $\sigma$ and $r>1$ the inequality
\begin{equation}\label{ineq:lemma2}
\sigma r > (b-\sigma)(b-1)
\end{equation}
is satisfied. Then the domain in which it is not satisfied is bounded on $\sigma$.

Condition \eqref{ineq:lemma2} is equivalent to
\begin{equation}\label{ineq:afterDivision}
  r > \frac{1}{\sigma}b^2 + \frac{-\sigma-1}{\sigma}b + 1.
\end{equation}
The right-hand side of the inequality is quadratic polynomial in $b$
with a positive coefficient of $b^2$. The abscissa of parabola peak is equal to
$b'_0 = \frac{\sigma}{2} + \frac{1}{2}$. If we consider $\sigma > 7$,
then $b'_0 > 4$.
In this case the maximum of the right-hand side of inequality
for  $b \in (0, 4)$
is smaller than the value for  $b = 0$, i.e. $1$. Thus,
for  sufficiently great $\sigma$, $b \in (0, 4)$ and $r>1$
the relation
\begin{equation*}
r > 1 > \frac{1}{\sigma}b^2 + \frac{-\sigma-1}{\sigma}b + 1
\end{equation*}
is satisfied, i.e. \eqref{ineq:lemma2} is valid.

\textbf{III.}
If $r > 1$, $b \in (0, 4)$, and $\sigma > 7$,
then condition \eqref{cond:sufficientLemma2}
and inequality \eqref{ineq:lemma2} are satisfied,
i.e. the conditions of Theorem \ref{Main} are satisfied.
\end{proof}

\begin{proof} (Remark 3)
 This is obvious from conditions \eqref{cond:sufficientLemma2}.
\end{proof}

\section*{Acknowledgements}
 This work was supported by the Russian Scientific Foundation (project 14-21-00041)
 and Saint-Petersburg State University.

\section*{References}
%\bibliographystyle{elsarticle-num}
%\bibliography{C:/Dropbox/bib/bib_full,C:/Dropbox/bib/bib_leonov,C:/Dropbox/bib/bib_nk}
%,C:/Dropbox/bib/genlorenz-bib,C:/Dropbox/bib/bib-2008-str-at

\end{document}